\documentclass[11pt,reqno]{article}

\usepackage[usenames]{color}
\usepackage{graphicx}
\usepackage{amscd}
\usepackage[english]{babel}
\usepackage{color}
\usepackage{fullpage}
\usepackage{psfig}
\usepackage{graphics,amsmath,amssymb}
\usepackage{amsthm}
\usepackage{amsfonts}
\usepackage{latexsym}
\usepackage{epsf}
\usepackage{float}
\usepackage{MnSymbol}
\usepackage{tikz} 
\usepackage{ifthen} 
\usepackage{lscape} 
\usepackage{verbatim} 
\usetikzlibrary{arrows,decorations.markings}

\usepackage{hyperref}

\definecolor{webgreen}{rgb}{0,.5,0}
\definecolor{webbrown}{rgb}{.6,0,0}
		
\newcommand{\R}{{\mathbb R}}

\theoremstyle{plain}
\numberwithin{equation}{section}

\newtheorem{remark}{Remark}[section]

 \newcommand{\seqnum}[1]{\href{http://oeis.org/#1}{\underline{#1}}}

\newcommand{\eqn}[1]{(\ref{#1})}

\newcommand{\ta}{\theta} 

\newcommand{\eps}{\varepsilon}
\newcommand{\beql}[1]{\begin{equation}\label{#1}}
\newcommand{\eeq}{\end{equation}}

\begin{document}

\setcounter{page}{0}

 \newcounter{VarX}
 \setcounter{VarX}{0}


\begin{center}

{\large\bf On Dissecting  Polygons into  Rectangles } \\
\vspace*{+.2in}
\
N.\ J.\  A.\  Sloane \\ 
The OEIS Foundation Inc., 
11 South Adelaide Ave.,
Highland Park, NJ 08904, USA \\
Email:  \href{mailto:njasloane@gmail.com}{\tt njasloane@gmail.com}

\vspace*{+.1in}

Gavin A.\ Theobald \\
15 Glasdrum Rd., Fort William, Inverness-shire, Scotland PH33 6DD, UK \\
 Email:  \href{gavintheobald@icloud.com }{\tt gavintheobald@icloud.com }

\vspace*{+.1in} 




\vspace*{+.1in}
\begin{abstract}
What is the smallest number of pieces that you can cut an  $n$-sided regular polygon into
 so that the pieces can be rearranged to form a rectangle? Call it $r(n)$.
The rectangle may have any proportions you wish, as long as it is a rectangle. 
The rules are the same as for 
 the classical problem where the rearranged pieces must form a square. Let $s(n)$ denote the minimum
 number of pieces for that problem.
For both problems the  pieces may be turned over and the cuts must be simple curves.
The conjectured values of $s(n), 3 \le n \le 12$, are $4, 1, 6, 5, 7, 5, 9, 7, 10, 6$.
However, only $s(4)=1$ is known for certain. 
The problem of finding $r(n)$ has received less attention.
In this paper we give constructions showing that $r(n)$ for $3 \le n \le 12 $ is at most 
$2, 1, 4, 3, 5, 4, 7, 4, 9, 5$, 
improving on the bounds for $s(n)$ in every case except $n=4$.
For the $10$-gon our construction uses three fewer pieces than the bound for $s(10)$.
Only $r(3)$ and $r(4)$ are known for certain. We also briefly discuss $q(n)$, 
the minimum number of pieces needed to dissect a regular $n$-gon into a monotile.
\end{abstract}
\end{center}

\section{Introduction}\label{Sec1}
Two polygons are  said to be {\em equidecomposable} if one can be cut into a finite number of pieces that can be rearranged to form the other.
Pieces may be turned over, and the cuts must be along simple plane curves.
The Bolyai-Gerwien theorem from the 1830s states that any two polygons of the same area are
equidecomposable, and the dissection can be carried out using only
triangular pieces. Furthermore, the dissection can be carried out using only a straightedge and compass.
Boltianskii \cite{Bol78} gives an excellent survey.

\begin{figure}[!ht]
\centerline{\includegraphics[angle=0, width=2.5in]{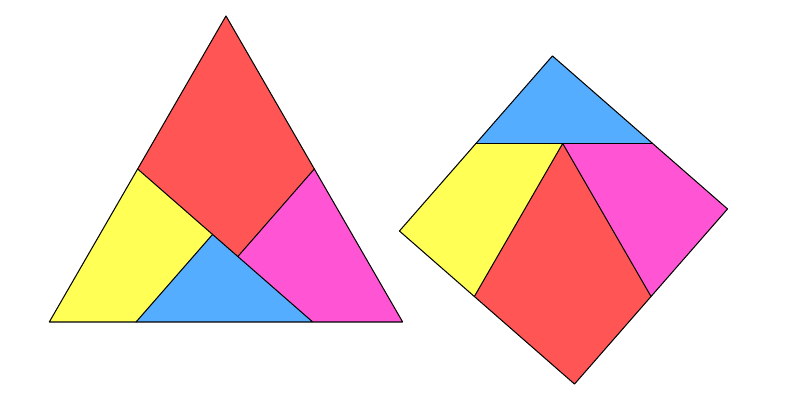}}
\caption{A $4$-piece dissection of an equilateral triangle into a square.}
\label{Fig3to4}
\end{figure}

A much-studied special case of this problem asks for
the minimum number of pieces ($s(n)$, say) of any shape
needed to dissect a regular polygon with $n$ sides into a square
of the same area.
Despite its long history
\cite{AkMa17, EppGJ, Fred97, Fred97a, Fred02, Fred06, Fred18, Lind64, Lind72, MRE39, GDDb}, surprisingly little is known about this problem.
The best upper bounds currently known for $s(n), n = 3, 4, 5, \ldots, 16$, are shown in Table~\ref{Tab1}.
The only exact value known appears to be the trivial result that $s(4) = 1$.
Even the value of $s(3)$ is not known. The four-piece dissection of an equilateral triangle into a square 
shown in Fig.~\ref{Fig3to4} is at least $120$ years old
 (see the discussions in \cite{AkMa17},  \cite[Ch.\ 12]{Fred97}), but there is no proof that it is optimal.
 The conjecture that it is impossible to dissect an equilateral triangle into three pieces 
that can be rearranged to form a square must
 be one of the oldest unsolved problems in geometry.

\begin{table}[htb]
$$
\begin{array}{|c|rrrrrrrrrrrrrr|}
\hline
n & 3 & 4 & 5 & 6 & 7 & 8 & 9 & 10 & 11 & 12 & 13 & 14 & 15 & 16 \\
\hline
s(n) \le           & 4 & 1 & 6 & 5 & 7 & 5 & 9 & 7 & 10 & 6  & 11 & 9 & 11 & 10 \\ \hline
r(n) \le            & 2 & 1 & 4 & 3 & 5 & 4 & 7 & 4 &  9  & 5 & 10 & 7 & 10 &   9 \\ \hline
q(n) \le           & 1 & 1 & 2 & 1 & 3 & 2 & 3 & 2 &  4  & 3 &   4 & 3 &  5  &   4 \\ \hline
\end{array}
$$
\caption{$s(n)$ (resp.\ $r(n)$, $q(n)$) is the  minimum number of pieces needed to   
dissect a regular $n$-sided polygon into a square (resp.\ rectangle, monotile).
Only the values of  $s(4)$, $r(3)$, $r(4)$ 
and $q(n)$ for $n=3, 4, 5, 6, 8, 10$ are known to be exact.}\label{Tab1}
\end{table}

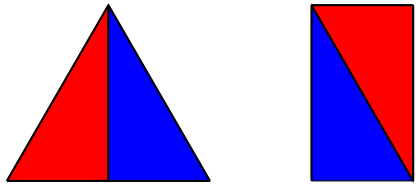
\begin{figure}[!ht]
\begin{center}
\begin{tikzpicture}[scale=2.7]
\def \SX{1};
\coordinate (A) at (-.5,0);
\coordinate (B) at (.5,0);
\coordinate (D) at (0,0);
\coordinate (C) at (0,0.866);
\coordinate (e) at (\SX + .5,.866);
\coordinate (b) at (\SX+.5,0);
\coordinate (d) at (\SX+0,0);
\coordinate (c) at (\SX+0,0.866);

\draw[fill = red] (A) -- (D) -- (C) -- (A);
\draw[fill = blue] (D) -- (B) -- (C) -- (D);
\draw[fill = blue] (d) -- (b) -- (c) -- (d);
\draw[fill = red] (c) -- (b) -- (e) -- (c);
\draw[thick] (b) -- (e);
\draw[thick] (A) -- (B);
\draw[thick] (A) -- (C);
\draw[thick] (D) -- (C);
\draw[thick] (B) -- (C);
\draw[thick] (d) -- (c);
\draw[thick] (b) -- (c);
\draw[thick] (b) -- (e);
\draw[thick] (c) -- (e);
\end{tikzpicture}
\caption{A $2$-piece dissection of an equilateral triangle into a rectangle. One piece must be turned over.}
\label{Fig3gon} 
 \end{center}
 \end{figure}

\begin{figure}[!ht]
\begin{center}
\begin{tikzpicture}[scale=1.5]
\def \SX{2.5};

\coordinate (A) at (-1,0);
\coordinate (B) at (-.5,0);
\coordinate (C) at (1,0);
\coordinate (D) at (-.5, 0.866);
\coordinate (G) at (-.5, -0.866);
\coordinate (E) at (.5, 0.866);
\coordinate (H) at (.5, -0.866);
\coordinate (F) at (1, .866);
\coordinate (I) at (1, -.866);

\coordinate (b) at (\SX+-.5,0);
\coordinate (c) at (\SX+1,0);
\coordinate (d) at (\SX-.5, 0.866);
\coordinate (g) at (\SX-.5, -0.866);
\coordinate (e) at (\SX+.5, 0.866);
\coordinate (h) at (\SX+.5, -0.866);
\coordinate (f) at (\SX+1, .866);
\coordinate (i) at (\SX+1, -.866);

\draw[fill = red] (G) -- (H) -- (C) -- (E) -- (D) -- (G);
\draw[fill = blue] (A) -- (B) -- (D) -- (A);
\draw[fill = yellow] (A) -- (B) -- (G) -- (A);
\draw[fill = red] (g) -- (h) -- (c) -- (e) -- (d) -- (g);
\draw[fill = blue] (A) -- (B) -- (D) -- (A);
\draw[fill = yellow] (A) -- (B) -- (G) -- (A);
\draw[fill = blue] (h) -- (i) -- (c) -- (h);

\draw[fill = yellow] (c) -- (e) -- (f) -- (c);

\draw[thick] (A) -- (B);
\draw[thick] (A) -- (D);
\draw[thick] (A) -- (G);
\draw[thick] (D) -- (G);
\draw[thick] (D) -- (E);
\draw[thick] (C) -- (E);
\draw[thick] (C) -- (H);
\draw[thick] (G) -- (H);

\draw[thick] (d) -- (f);
\draw[thick] (f) -- (i);
\draw[thick] (i) -- (g);
\draw[thick] (g) -- (g);

\end{tikzpicture}
\caption{A $3$-piece dissection of a regular hexagon into a rectangle.}
\label{Fig6gon} 
\end{center}

\end{figure}
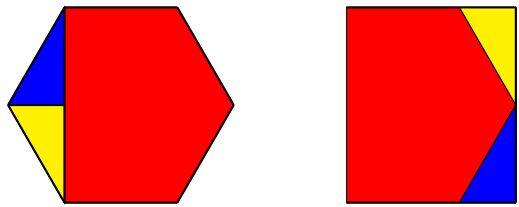

In the present article we consider a weaker constraint: what
is the minimum number of pieces ($r(n)$, say) 
needed to dissect a regular polygon with $n$ sides into a {\em rectangle}
of the same area?  The proportions of the rectangle can be anything you want.
It is clear that $r(3)=2$ and $r(6) \le 3$ (Figs.~\ref{Fig3gon}, \ref{Fig6gon}).  
(Surely it should be possible to prove 
that no two-piece dissection of a regular hexagon into a rectangle is possible?)

Lindgren \cite{Lind64, Lind72} gives examples of regular $n$-gon to non-square rectangle dissections,
but none have fewer pieces than the corresponding $n$-gon to square dissections.
Frederickson~\cite[pp.\ 150-151]{Fred97} mentions that in 1926 H.~E.~Dudeney found a $4$-piece octagon to rectangle dissection (probably that shown in Fig.~\ref{Fig8A}).

In April 2023, before the current investigation was begun, the second author (G.A.T.)'s
{\em Geometric Dissections} database~\cite{GDDb} contained 
 several examples of regular $n$-gon to rectangle dissections that 
had fewer pieces than the best square dissections known, including
a $5$-piece dissection of a pentagon, a $4$-piece dissection of an octagon,
a $6$-piece dissection of a $10$-gon, and a $5$-piece dissection of a $12$-gon.
The database also contained the star polygon 
and  Greek cross dissections shown below  in \S\ref{SecStar} 
and \S\ref{SecGC}.\footnote{\cite{GDDb} now includes all the dissections mentioned in this article.} 

In May 2023, Adam Gsellman \cite{Gsel23} wrote to N.J.A.S.\  enclosing dissections
showing that $r(5) \le 4$, $r(7) \le 6$, and $r(8) \le 4$. 
As mentioned above, the third of these results was already known, 
and G.A.T.\ was aware that $r(5) \le 4$, although that result 
had not yet been mentioned in the literature, but the upper bound on $r(7)$ appeared to be new.
We have since shown that $r(7) \le 5$ (see~\S\ref{Sec7}),
but Gsellman's dissections of the pentagon and octagon are shown 
in~\S\ref{Sec5} and \S\ref{Sec8}. 

Table~\ref{Tab1} shows the current best upper bounds on $s(n)$ and $r(n)$ for $n \le 16$,
although to avoid making the paper too long we shall say very little about the cases when
$n > 12$.

The recent remarkable 
 discovery \cite{SMKG23} of a single polygon, or {\em monotile}, that tiles the plane,
but can only do so in a  non-periodic way, reminded us of a question asked by
Gr\"{u}nbaum and Shephard in 1986 \cite[\S2.6]{GS86}: 
what is the minimum number of pieces needed to dissect a
regular $n$-gon into a monotile that tiles the plane (allowing periodic tilings)?
Call this number $q(n)$.  Of course
squares and rectangles themselves are monotiles, so
$q(n) \le r(n) \le s(n)$. 
We have included  bounds on $q(n)$ in Table~\ref{Tab1}.
All  of $s(n)$, $r(n)$, and $q(n)$ are fundamental geometric quantities, with $q(n)$
perhaps the most basic of the three.
We will give some examples of these monotiles in Remark~\ref{RemMonotile},
including a dissection of a $9$-gon which
improves on Gr\"{u}nbaum and Shephard's. 

The paper is arranged as follows.
The remainder of this section gives some general remarks about our dissections.
Section~\ref{Sec2} defines some parameters and coordinates that will be generally useful 
for our dissections.
Subsequent sections deal in turn with pentagons, heptagons, octagons, up through $16$-gons, followed by
 sections  on star polygons and the Greek cross. Finally, 
Section~\ref{SecCurved} gives some examples of dissections 
where curved cuts appear essential if one wishes to minimize the number of pieces.
Dissections given without attribution  are believed to be new.

\begin{remark}\label{RemPf}
Certifying the dissections.
\end{remark}
\noindent We have attempted to give detailed descriptions of most of the dissections in the
main body of the paper (\S\ref{Sec5}-\S\ref{Sec14}), enough at any rate to convince
the reader that the dissections are correct. If the dissection begins
by cutting up  the regular polygon,
for example,  we have to make sure that the rearranged pieces 
form a proper rectangle. The pieces must not overlap, there can be no holes;
when pieces fit together at a vertex, the sum of the angles
must be $2\pi$ at an interior point, or $\pi$ or $\pi/2$ at a boundary point, and so on.

In simple cases the the correctness can be checked ``by hand'', 
like solving a jigsaw puzzle. The first pentagon dissection, in \S\ref{Sec5A},
is an example.

Many of our dissections were obtained by one of the standard strip or superposition constructions.
There are a great many versions of these constructions, and they are described
in most of the books on the subject,
and in the {\em Methods} section of \cite{GDDb},
 so we shall not say much about them here.\footnote{For
the mathematical theory underlying these constructions, see \cite{AFF00, Mul88, Schobi}.}

A simple example of a superposition, used to dissect a polygon $A$ into a polygon $B$,
 is to overlay a strip tiled with pieces from  $A$ and a second strip tiled with pieces
 from $B$. Then with luck the intersection of the two strips
 will provide the desired dissection  (see for example \cite[Ch.~11]{Fred97}, \cite[Chaps. 2-5]{Lind64}.)
 Since for our problem we can assume $B$ is a rectangle, we can often dispense with the second strip.
 All we need then  is a strip tiled with pieces from $A$. We obtain the dissection
 by cutting out a rectangle of the correct length from the strip.
 Examples are shown in Figs.~\ref{Fig7gonStripB} and \ref{Fig9B}.

If a dissection is obtained by one of these standard constructions, it can generally
be assumed to be correct. However, one must be careful: with complicated strips like those shown
in sections \S\ref{Sec9} onwards, it is easy to be mistaken about points
coinciding, or when triangular regions shrink to a point (in order to save a piece). 

For this and other reasons, we have  therefore tried to give {\em ab initio} descriptions of 
the dissections.  In most cases we are able to give  a straightedge and compass construction,
and to sketch a proof that it is correct.

\begin{remark}\label{RemSEC}
Straightedge and compass constructions.
\end{remark}
\noindent 
Given an initial drawing of a regular $n$-gon, our dissections can  be constructed using only a straightedge and compass.
That is, there is no need for a ruler: the construction does not require locating a point 
which is at some arbitrary irrational  distance from another point.\footnote{Hadlock~\cite{Had78} contains an excellent introduction to straightedge and compass constructions.}
We need to be given the initial $n$-gon, since, for example, a regular heptagon cannot be 
obtained with only a straightedge and compass.

Besides its aesthetic appeal, the advantage of a straightedge and compass construction
is that it enables us to give explicit coordinates 
for every vertex in the construction.\footnote{Although the Bolyai-Gerwien theorem 
guarantees that a straightedge and compass dissection of a regular
$n$-gon to a square exists, we don't know that this is 
true for a dissection with the minimum number of pieces.}
We usually start from from the vertices \eqn{Pcoords}  of the $n$-gon, and every subsequent 
vertex is then determined. In \S\ref{Sec10} we start from the rectangle, which we assume has 
width $\sqrt{5}$ and height $\cos(\pi/10)$, and again all subsequent vertices are determined.

Although in theory we could find exact expressions for all the coordinates
in a straightedge and compass construction  in this way, the expressions rapidly become unwieldy. 
In practice we have found it better to use computer algebra
systems such as WolframAlpha and Maple to guess expressions for the coordinates 
based on $20$-digit decimal expansions,
knowing that they could be justified if necessary.
We shall see examples of this is \S\ref{Sec9}.

\begin{remark}\label{RemTurn}
Turning pieces over.
\end{remark}
\noindent Although the definitions of $s(n)$ and $r(n)$ allow pieces to be turned over, this is deprecated by purists.
Fortunately all the dissections $s(n)$ and $r(n)$ for  $3 \le n \le 16$ mentioned
in Table~\ref{Tab1} can be accomplished without turning
pieces over, with the single exception of $r(3)$ (see Fig.~\ref{Fig3gon}),  which seems to require three pieces if turning over is forbidden.

Another example where turning pieces over appears to be essential to achieve the minimum 
number of pieces is the seven-piece dissection of $\{6\}$ into $\{8\}$ given in~\cite{GDDb}.

\begin{figure}[!ht]
\centerline{\includegraphics[angle=180, width=3.0in]{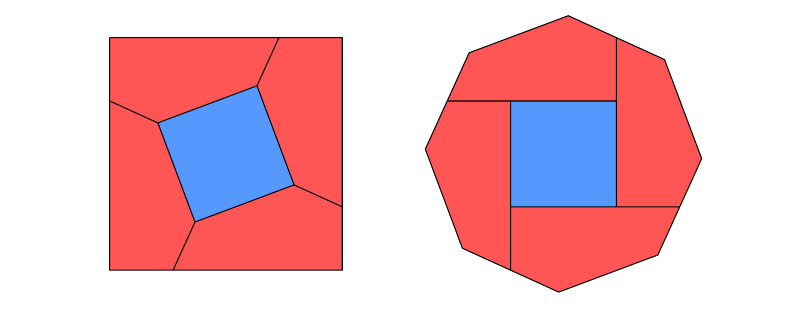}}
\caption{A $14$-th century $5$-piece regular octagon to square dissection, with cyclic $4$-fold symmetry
in both the octagon and the square.}
\label{FigOctSquare}
\end{figure}

\begin{remark}\label{RemCvx}
Convex pieces.
\end{remark}
\noindent 
The  dissections in Figs.~\ref{Fig3to4}-\ref{FigOctSquare},
\ref{Fig5gonF}, \ref{Fig5gon2},
\ref{Fig8A}, \ref{Fig12L} and \ref{FigPentagram}
use only convex pieces.  
Other things being equal, we prefer convex pieces, of roughly equal size, that do not need to be turned over.
The primary goal however is always  to minimize the number of pieces.

\begin{figure}[!ht]
\centerline{\includegraphics[angle=180, width=3.0in]{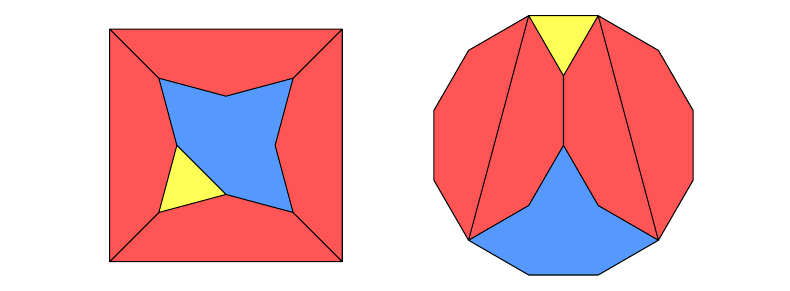}}
\caption{Lindgren's astonishing  $6$-piece $12$-gon to square dissection, with mirror symmetry \cite[Ch.\ 9]{Lind64}.}
\label{Fig12to4}
\end{figure}

\begin{remark}\label{RemClassical}
Improving on classic dissections.
\end{remark}
\noindent We were surprised to find that $r(8)$ is apparently less than $s(8)$, and 
$r(12)$ apparently less than $s(12)$,
since the best octagon to square\footnote{The octagon in the well-known Chase Bank logo is different from the octagon in Fig.~\ref{FigOctSquare}. The Chase octagon has a square surrounded by four trapezoids, whereas  Fig.~\ref{FigOctSquare} has a square surrounded by four pentagons. The  pieces in the Chase logo can be rearranged to form a rectangle, but not a square.}  
and $12$-gon to square constructions are so
striking (Figs.~\ref{FigOctSquare},~\ref{Fig12to4}). 
One feels that they could not possibly be improved on.
Yet if we only want a rectangle, there is a four-piece dissection of the octagon that has
essentially the same symmetry as Fig.~\ref{FigOctSquare},
as we shall see in~\S\ref{Sec8}.
Likewise, for the $12$-gon, we can save a piece if we only want a rectangle, at the cost however of giving up 
all symmetry (see~\S\ref{Sec12}).

\begin{figure}[!ht]
\centerline{\includegraphics[clip=true, trim={4cm, 0cm, 4cm, 0},  angle=90,  width=0.4\linewidth]
{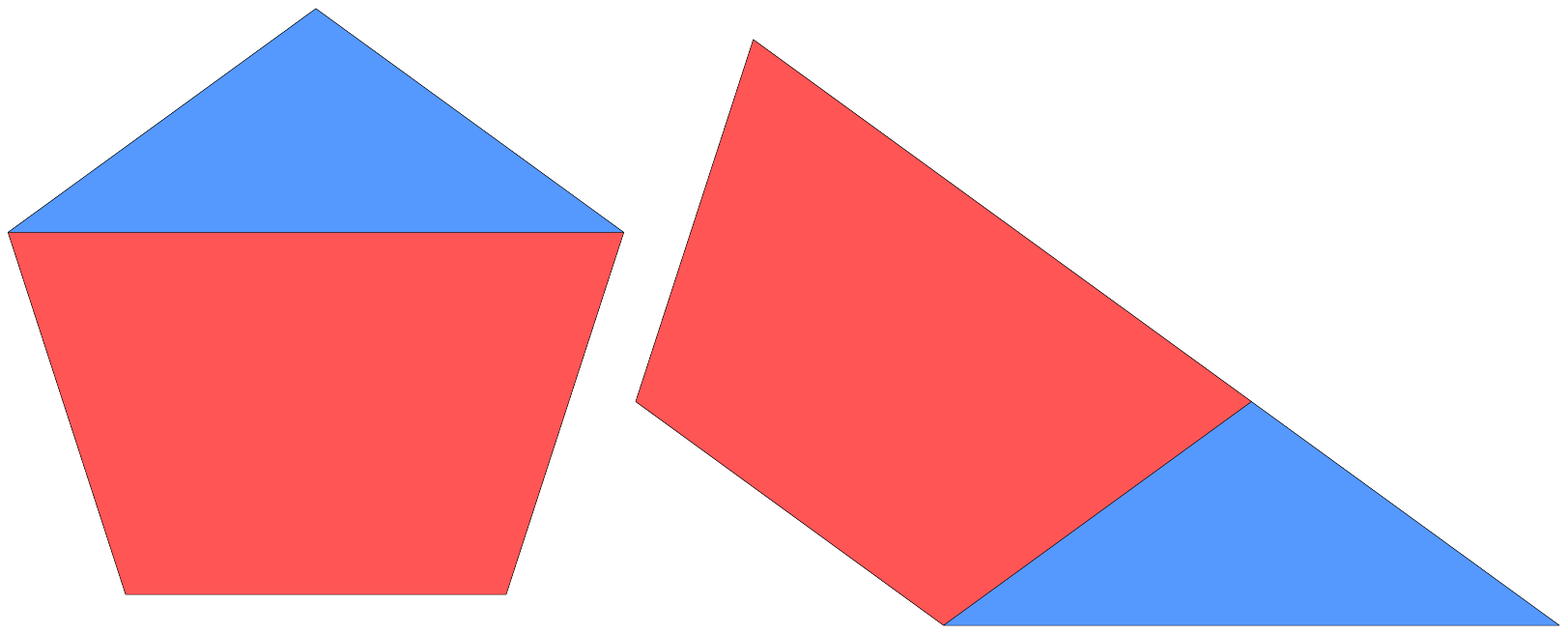}}
\caption{A two-piece dissection of a pentagon into a monotile, illustrating $q(5)=2$.} 
\label{Fig5gonMono}
\end{figure} 

\begin{figure}[!ht]
\centerline{\includegraphics[clip=true, trim={3cm, 0cm, 3cm, 0},  angle=90,  width=0.4\linewidth]
{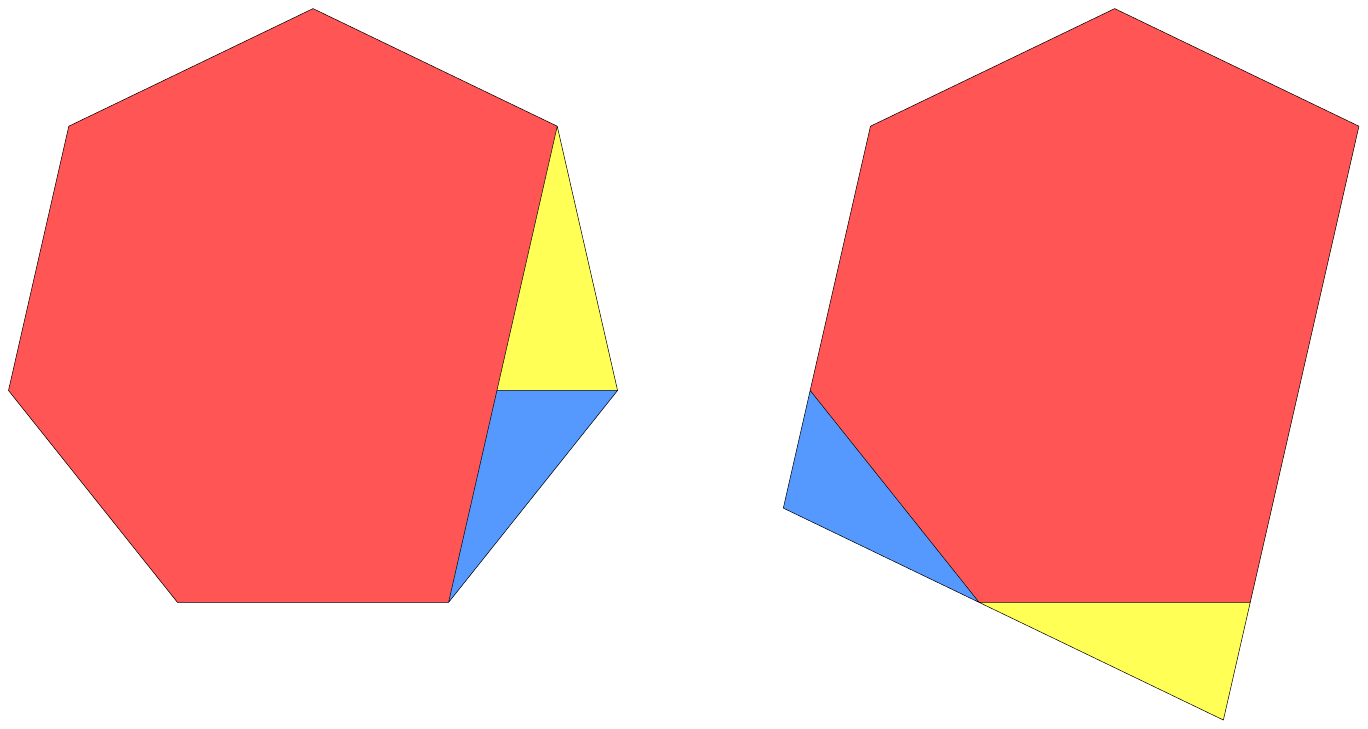}}
\caption{A three-piece dissection of a heptagon into a monotile, illustrating $q(7) \le 3$ \cite{Lind64}.} 
\label{Fig7gonMono}
\end{figure}

\begin{figure}[!ht]
\centerline{\includegraphics[clip=true, trim={2cm, 0cm, 2cm, 0},  angle=90,  width=0.42\linewidth]
{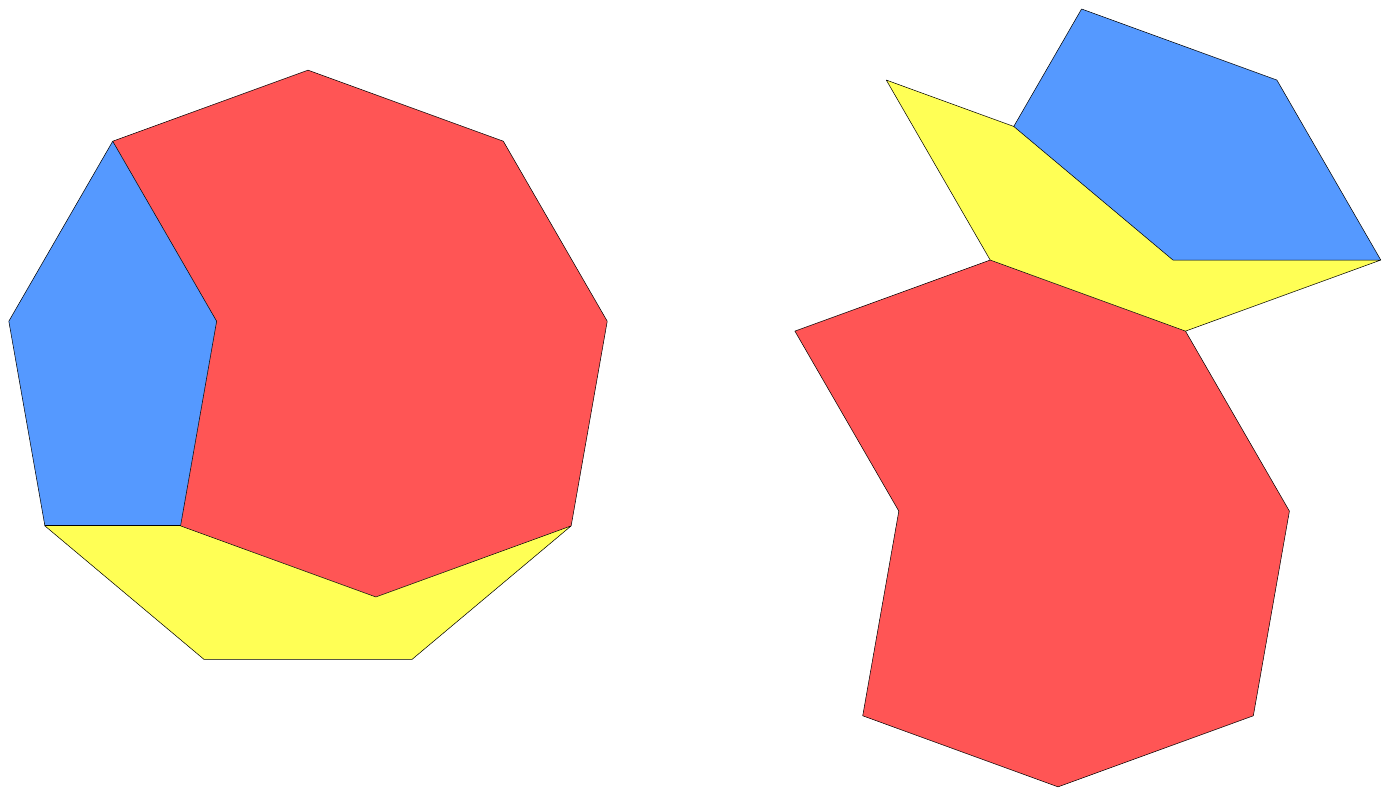}}
\caption{A three-piece dissection of a $9$-gon into a monotile, illustrating $q(9) \le 3$. This improves on the   dissection in \cite[Fig. 2.6.1]{GS86}.} 
\label{Fig9gonMono}
\end{figure}

\begin{figure}[!ht]
\centerline{\includegraphics[clip=true, trim={2cm, 0cm, 2cm, 0},  angle=90,  width=0.4\linewidth]
{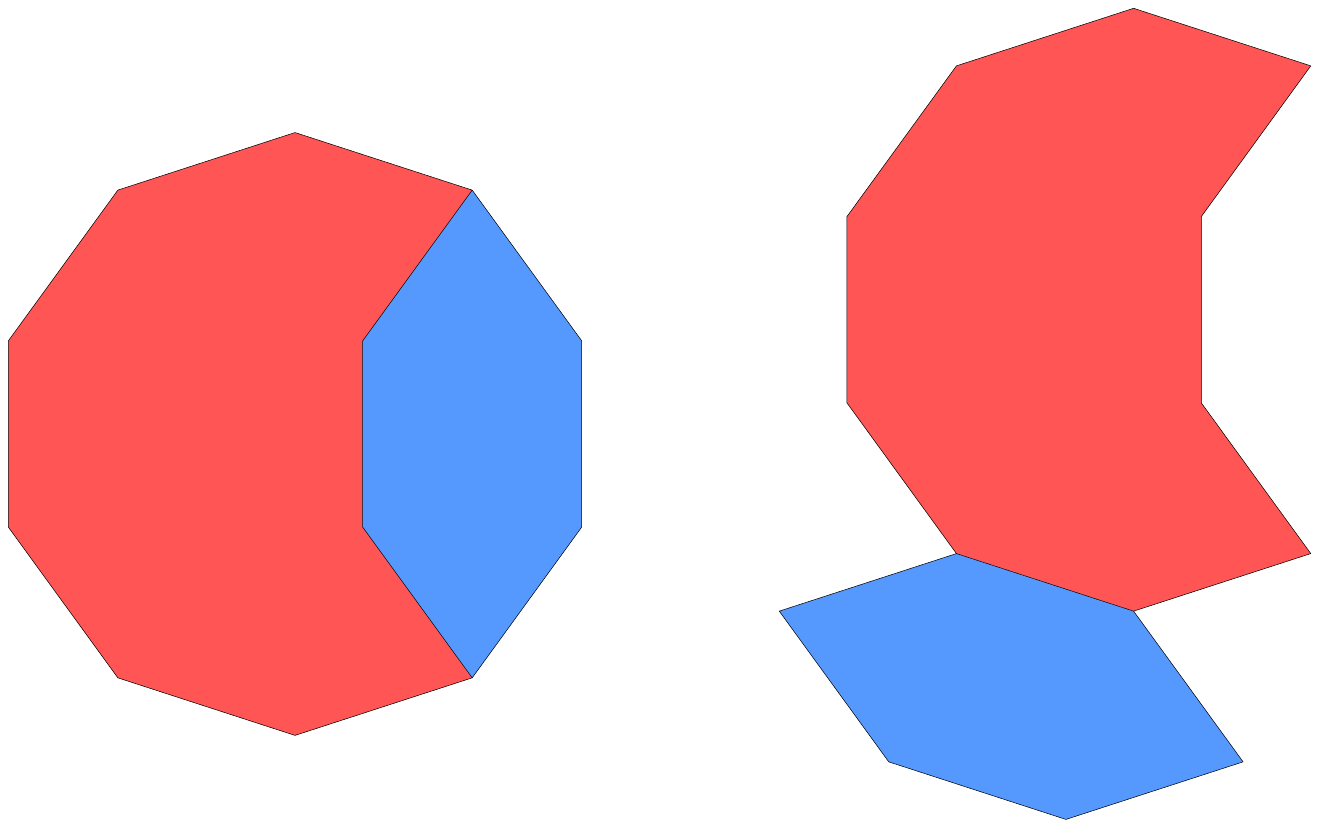}}
\caption{A two-piece dissection of a $10$-gon into a monotile, illustrating $q(10) = 2$.} 
\label{Fig10gonMono}
\end{figure}

\begin{remark}\label{RemMonotile}
Dissecting a regular $n$-gon into a monotile.
\end{remark}
\noindent As mentioned above, the problem of finding $q(n)$ is discussed in Section~2.6 of
Gr\"{u}nbaum and Shephard  \cite{GS86},
where a dissection attaining $q(n)$ is called a {\em minimal dissection tiling}.
Figure 2.6.1 of \cite{GS86} shows 
dissections achieving $q(5) = q(8) = q(10) = 2$
and conjectured solutions for $n = 7, 9, 12$, with a reference to Lindgren \cite{Lind72}.
Our Figures~\ref{Fig5gonMono}, \ref{Fig7gonMono}, \ref{Fig9gonMono}, and
\ref{Fig10gonMono} show examples of monotiles for $n = 5, 7, 9$, and $10$.
Many further examples with $n \le 17$  can be seen in OEIS entry \seqnum{A362938}.

By generalizing the constructions for $n = 8, 10, 12, 14, \ldots$, we can
show that $q(2t) \le \lfloor t/2\rfloor$ for $t \ge 2$, and it may even be true
that  $q(2t) = \lfloor t/2\rfloor$ holds for all $t \ge 2$. No similar conjecture is known
for $q(2t+1)$.

\begin{remark}\label{RemLB}
Lower bounds.
\end{remark} 
\noindent It seems to be difficult to obtain 
nontrivial lower bounds 
on  any of $s(n)$, $r(n)$, or $q(n)$.\footnote{G.A.T. conjectures that
$1/4 \le q(n)/n \le r(n)/n \le s(n)/n$ for $n>10$, and that all three of these ratios approach $1/4$ as $n$ increases.}
References \cite{CKU07, KKU98}  do give some lower bounds, but only 
for more restricted classes of dissections (only allowing polygonal cuts, or
alternatively  what are called ``glass cuts'').
Two  negative results 
we do know of are the easily-proved fact\footnote{An equilateral triangle
of area $1$ has side-length $2/3^{1/4} = 1.519$,
which won't fit into a square of area $1$. So each vertex of the triangle must be in a different piece of the dissection.}  that $s(3)$ cannot be $2$, and 
the nontrivial theorem  \cite{DHK63} that a circular disk 
cannot be dissected into a polygon.

We cannot resist mentioning that the latter question has been in the news recently, because of  progress on 
the problem of dissecting a circular disk into a square
{\em if fractal pieces are allowed} \cite{MNP22}. The new construction involves at most $10^{200}$ pieces.

\begin{remark}\label{RemOEIS}
The OEIS entry.
\end{remark}
\noindent The best upper bounds currently known for 
$s(n)$, $r(n)$, and $q(n)$  for $n$ up to around $16$ or $20$
are  listed 
in the {\em On-Line Encyclopedia of Integer Sequences} (or {\em OEIS}) database~\cite{OEIS}
as sequences \seqnum{A110312},  \seqnum{A362939}, and \seqnum{A362938}, respectively.
This is an exception to the usual OEIS policy of requiring that all terms in sequences must be known 
exactly, but these sequence are included because of their importance and in the hope 
that someone will establish the truth of some of the conjectured values.

\begin{remark}\label{RemApp}
Applications
\end{remark}
\noindent These polygon to rectangle dissections have  potential applications
to lossless source coding (cf.~\cite{Schobi}).
If a source (a lens, perhaps) repeatedly produces an output which is a point 
in a regular $12$-gon, say, then the dissection  could be used to map the point to a
more convenient  pair
$(x,y)$ of rectangular coordinates (compare Figs~\ref{Fig12L} and \ref{Fig12B}).

\begin{remark}\label{RemFurther}
For further information.
\end{remark}
\noindent
The database~\cite{GDDb} is the best reference for drawings of
dissections mentioned in Table~\ref{Tab1} but not included in the present article.

\begin{remark}\label{RemNotation}
Notation
\end{remark}
\noindent For up to eight sides we will use the names  triangle, $\ldots$, hexagon, heptagon, and octagon, but for nine or more sides
we will say $9$-gon, $10$-gon, $\ldots$. 
The symbol $\{n\}$  denotes a regular $n$-sided polygon, and 
 $\{n/m\}$ is a regular star polygon (cf.~\cite{Cox61}).
 $L_{k,n}$ is the length of the chord joining two vertices of $\{n\}$ that are $k$ edges apart  \eqn{Lchords}.
 Decimal expansions have been truncated not rounded.

\section{Regular polygons: coordinates and metric properties}\label{Sec2}

In later sections we will usually begin with a regular $n$-sided polygon,  with  $n \ge 3$,
having sides of length $1$, with the goal of cutting it into as few pieces as possible 
which can  be rearranged to form a rectangle.

We often take the polygon to have center $P_0 = (0,0)$ at the origin:, 
and to have vertices $P_1, P_2, \ldots$, $P_n$ labeled counterclockwise, starting 
with $P_1$ at the apex of the figure (Fig.~\ref{FigNgon1}).
The angle subtended at the center by an edge is $2 \pi / n$, and 
we define $\theta = \pi/n$ and $\phi = \pi/2 - \theta = \frac{n-2}{2n}\pi$,
which will be the principal angles used in the formulas.

\ifnum \value{VarX}<1 


{
\begin{figure}[!ht]
\centerline{\includegraphics[angle=0, width=4in]{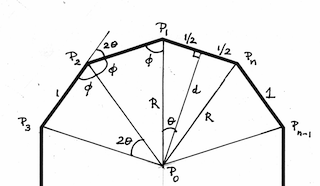}}
\caption{Vertices and principal angles and distances for a regular $n$-gon.}
\label{FigNgon1}
\end{figure} 
}
\else 
{
\begin{figure}[!ht]
\centerline{\includegraphics[angle=0, width=4in]{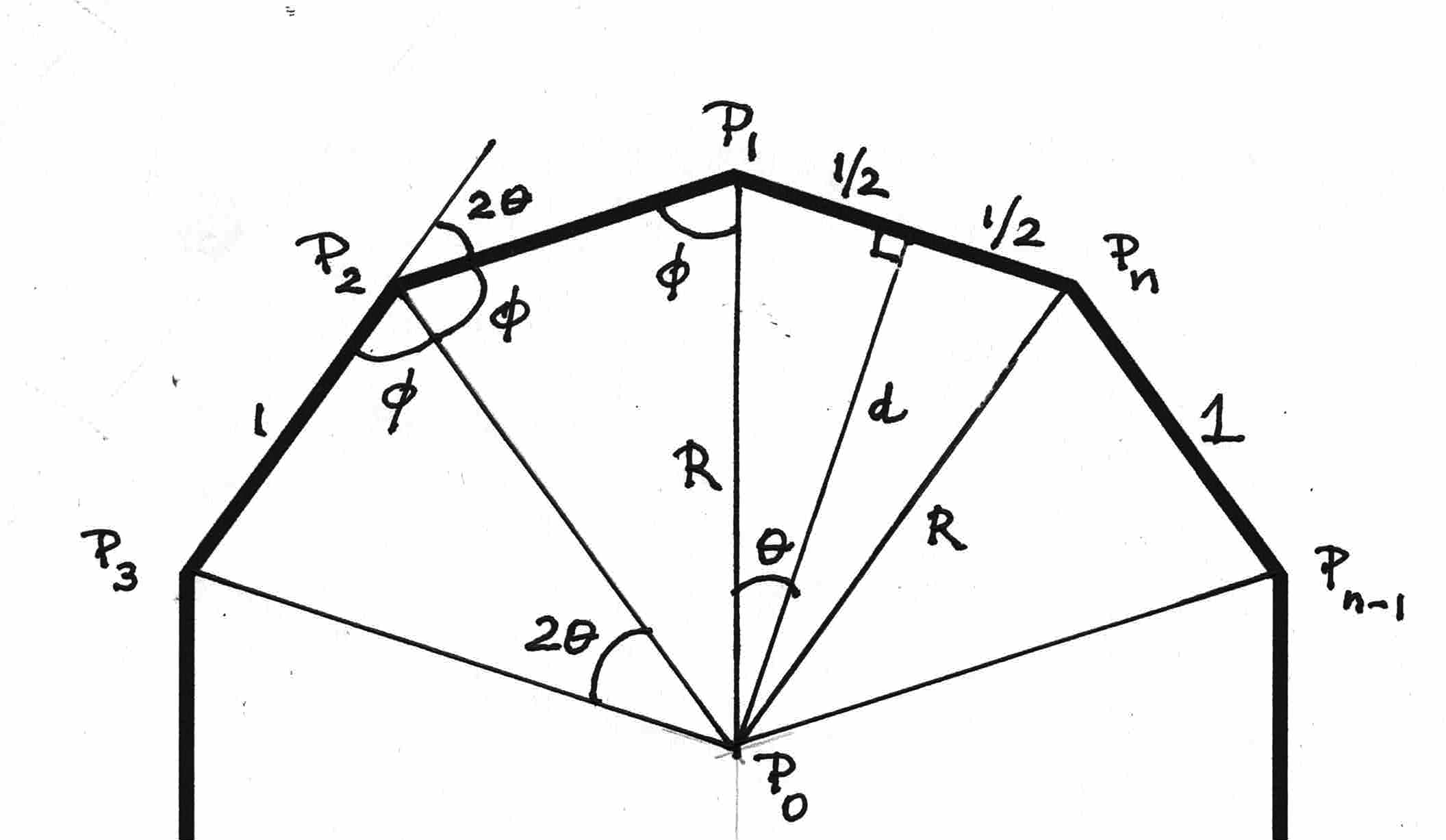}}
\caption{Vertices and principal angles and distances for a regular $n$-gon.}
\label{FigNgon1}
\end{figure} 
}
\fi

Since the sides have length $1$, the radius of the polygon is $R = \frac{1}{2 \sin \theta}$,
and the distance from the center to the midpoint of an edge is $d = \frac{1}{2 \tan \theta}$.
The vertices have coordinates
\begin{align}\label{Pcoords}
P_1  & = (0,R)\,, \nonumber \\
P_k  & = (-R \sin (2(k-1)\theta), R \cos (2(k-1)\theta))\,, \nonumber \\
P_{n+2-k}  & = (R \sin (2(k-1)\theta), R \cos (2(k-1)\theta))\,, \nonumber \\
\end{align}
for $k = 2, \ldots, \lfloor(n+2)/2 \rfloor$.
The polygon has area 
\beql{EqArea}
\frac{nd}{2} ~=~  \frac{n}{4} \cot \frac{\pi}{n}.
\eeq

Many dissections involve a cut along a chord $P_i P_{i+k}$ joining two vertices.
Assuming $n \ge 3$, 
$k \ge 0$, let  $L_k = L_{k,n}$ denote the length of the chord joining two vertices that are $k$
edges apart in $\{n\}$ and are in the same semicircle.
Then $L_0 = 0, L_1 = 1$,  $L_k = L_{k-2} + 2 \cos (k-1) \theta$, 
and for $t \ge 1$,
\beql{Lchords}
L_{2t}   = 2 \sum_{j = 0}^{t-1}  \cos (2j+1) \theta \,, \quad 
L_{2t+1}   = 1 + 2 \sum_{j = 1}^{t}  \cos 2j \theta \,.
\eeq 

\section{Four-piece dissections of a pentagon}\label{Sec5}

In 1891 Robert Brodie discovered a six-piece dissection of a regular pentagon into a square
(\cite[p.\ 120]{Fred97}, \cite[Fig.\ 3.1]{Lind64},
\cite{GDDb}), and there has been no improvement since then,
so it seems likely that $s(5)=6$.  In this section we give four different four-piece dissections
of a regular pentagon into a rectangle, showing that $r(5) \le 4$.
Rather surprisingly, this result does not seem to have been
mentioned in the literature before now.
The fact that there are at least  four ways to get $r(5) \le 4$ makes us wonder if $r(5)$ might
 actually be equal to $3$.

\begin{figure}[!ht]
\begin{center}
\begin{tikzpicture}[scale=3]
\def \R{.851};
\coordinate (P0) at (0,0);
\coordinate (P1) at (0,\R);
\coordinate (P2) at (-.809,.263);
\coordinate (P3) at (-.5,-.688);
\coordinate (P4) at (.5,-.688);
\coordinate (P5) at (.809,.263);
\coordinate (Q1) at (0,-.688);
\coordinate (A) at (.155,.740);
\coordinate (B) at (-.5,.263);
\coordinate (C) at (0,.263);
\coordinate (Cr) at (0.15,.263);
\coordinate (D) at (1.0,.263);
\coordinate (E) at (1.309,.263);
\coordinate (F) at (1.150,-.213);
\coordinate (G) at (1.309,-.688);
\draw[ultra thick] (P1) -- (P2);
\draw[ultra thick] (P2) -- (P3);
\draw[ultra thick] (P3) -- (P4);
\draw[ultra thick] (P4) -- (P5);
\draw[ultra thick] (P5) -- (P1);
\draw[ultra thick] (P2) -- (C);
\draw[ultra thick] (P1) -- (C);
\draw[ultra thick] (P2) -- (E);
\draw[ultra thick] (E) -- (G);
\draw[ultra thick] (G) -- (P3);
\draw[ultra thick] (B) -- (P3);
\draw[ultra thick] (D) -- (G);
\draw[ultra thick] (P4) -- (F);
\draw[ultra thick] (C) -- (A);
\draw[fill=red] (P1) -- (P2) -- (C) -- (A) -- (P1);
\draw[fill=red] (P4) -- (P5) -- (D) -- (F) -- (P4);
\draw[fill=yellow] (P3) -- (P2) -- (B) -- (P3);
\draw[fill=yellow] (G) -- (D) -- (E) -- (G);
\draw[fill=cyan] (C) -- (A) -- (P5) -- (C);
\draw[fill=cyan] (F) -- (G) -- (P4) -- (F);
\draw[fill=purple] (B) -- (P5) -- (P4) -- (P3) -- (B);
\draw[ thick] (P1) -- (Q1);
\node[above]  at  (P1) {$\mathbf{P_1}$};
\node[left]  at  (P2) {$\mathbf{P_2}$};
\node[below]  at  (P3) {$\mathbf{P_3}$};
\node[below]  at  (P4) {$\mathbf{P_4}$};
\node[above]  at  (P5) {$\mathbf{P_5}$};
\node[below]  at  (Q1) {$\mathbf{Q_1}$};
\node[right]  at  (A) {$~\mathbf{Q_3}$};
\node[above]  at  (B) {$\mathbf{Q_4}$};
\node[below]  at  (Cr) {$\mathbf{Q_2}$};
\node[above]  at  (D) {$~~\mathbf{Q_5}$};
\node[above]  at  (E) {$\mathbf{Q_8}$};
\node[left]  at  (F) {$\mathbf{Q_6}~$};
\node[below]  at  (G) {$\mathbf{Q_7}$};
\end{tikzpicture}
\caption{The pentagon $P_1 P_2 \ldots P_5$ is
cut into four convex pieces which can be rearranged to form the rectangle $Q_4 P_3 Q_7 Q_8$.}
\label{Fig5gonF}
\end{center}
\end{figure}
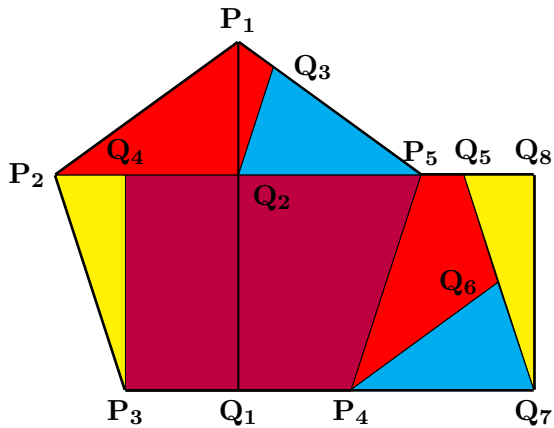

 The first two four-piece dissections (\S\ref{Sec5A}, \S\ref{Sec5B}) were 
 found by the authors (although they can hardly be new), 
and have the property that the pieces are convex; 
the other two (\S\ref{Sec5C}, \S\ref{Sec5D}) are due to Adam Gsellman \cite{Gsel23}. 

We use the notation introduced in the previous section
 (taking $n=5$).
 The pentagon has vertices $P_1, \ldots, P_5$
with sides of length $1$.
The key angles are $\ta = \angle P_1 P_2 P_5  = 36^\circ$, $\phi = 54^\circ$, 
$\angle P_2 P_1 P_5 = 2 \phi$, and
$\angle P_2 P_3 Q_4 = \theta/2$.
We note the values of
\begin{align}\label{EqPentagon1}
\sin 36^\circ & ~=~ \cos 54^\circ ~=~ \sqrt{ \frac{5-\sqrt{5}}{8}} \,, 
& \cos 36^\circ  ~=~ \sin 54^\circ ~=~ \frac{\sqrt{5}+1}{4}\,, \nonumber \\
\sin 72^\circ & ~=~ {\cos {18}}^{\circ}  ~=~   \sqrt{ \frac{5+\sqrt{5}}{8}} \,,
& \cos 72^\circ  ~=~ \sin 18^\circ ~=~ \frac{\sqrt{5}-1}{4}\,. 
\end{align}

\subsection{Pentagon \#1}\label{Sec5A}
To construct this dissection (see Fig.~\ref{Fig5gonF})  we draw a perpendicular from $P_1$ to 
the mid-point $Q_1$ of the opposite side, and draw the chord $P_2-P_5$.
Let $Q_2$ be the intersection of these two lines, and place $Q_3$ on $P_1-P_5$ so that 
$Q_2 P_5 Q_3$ is an isosceles triangle. Finally, draw a perpendicular $P_3-Q_4$
from $P_3$ to  $P_2-P_5$.

To form the rectangle, we first rotate the quadrilateral $P_1 P_2 Q_2 Q_3$  by $36^\circ$ 
and move it to $P_5 P_4 Q_6 Q_5$, then 
the triangle $Q_2 P_5 Q_3$  is moved to $Q_6 P_4 Q_7$, 
and the triangle $P_2 P_3 Q_4$ is moved horizontally to $Q_5 Q_7 Q_8$.
The two isosceles triangles (blue) have long sides of length $L_{2,5}/2 = (1+ \sqrt{5})/4$ and base $1/2$.
The two right triangles (yellow) have sides of lengths $\sin \theta$, $\cos \theta$, and $1$.

To  prove that this dissection is correct, we must check that, after the
 rearrangement, 
the result is indeed a rectangle, with area equal to that of the pentagon.
In particular, we must check that the points $P_2, Q_4, Q_2, P_5, Q_5, Q_8$ are collinear,
as are $P_3, P_4, Q_7$, and $Q_5, Q_6, Q_7$, that $Q_6$ bisects $Q_5-Q_7$, and also that the difference between the $x$-coordinates of 
$P_2$ and $P_3$ is equal to the difference
between the $x$-coordinates of $Q_5$ and $Q_8$.
Since we have complete information about the points, these checks are easily carried out.
The pentagon has area
$(5^{3/4}/(8  \sqrt{2}) ) (\sqrt{5}+1)^{3/2}$, and we can check that this is 
equal to $|P_3 Q_4|\cdot |P_3 Q_7|$. This completes the proof of the dissection.

It is worth pointing out that the trapezoid $P_5 P_4  Q_7 Q_5$ is symmetric about its vertical axis.

We chose this relatively simple example to illustrate the steps needed to prove that a
dissection is correct.  In later examples we will  just give the basic information needed for
the proof and leave the detailed verification to the reader.


\subsection{Pentagon \#2}\label{Sec5B}
This is very similar to the first dissection. Two of the pieces are the same, only now the rest of the pentagon is divided into two pieces that are reflections of each other. The trapezoid from Fig.~\ref{Fig5gonF}
has moved to the center of the rectangle.


\ifnum \value{VarX}<1


{
\begin{figure}[!ht]
\centerline{\includegraphics[angle=0, width=4in]{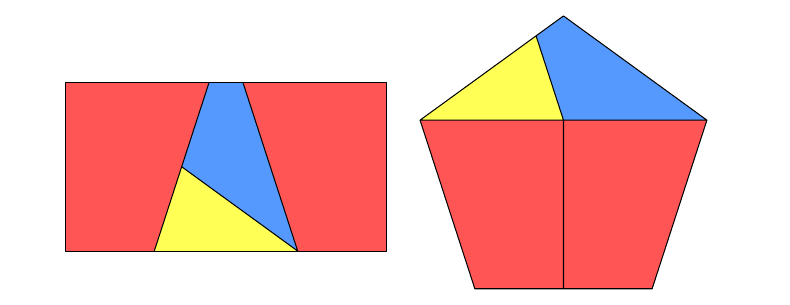}} 
\caption{Similar to the first pentagon dissection, only now the pieces are more nearly equal in size.}
\label{Fig5gon2}
\end{figure}
}
\else
{
\begin{figure}[!ht]
\centerline{\includegraphics[angle=0, width=4in]{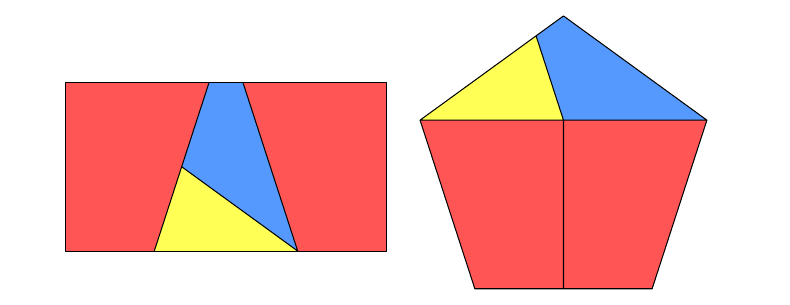}} 
\caption{Similar to the first pentagon dissection, only now the pieces are more nearly equal in size.}
\label{Fig5gon2}
\end{figure}
}
\fi

Although logically we are dissecting the polygon {\em into} a rectangle,
many of our colored illustrations have the rectangle on the left of the picture,
as in Figs.~\ref{Fig5gon2}, \ref{Fig7gon.d}, \ref{Fig8A}, etc.
This is because of the convention in  \cite{GDDb} of starting with
 the figure having the smaller number of vertices.
 Also, in the case of  strip or tessellation constructions, one often proceeds from the rectangle to the polygon.

\begin{figure}[!ht]
\centerline{\includegraphics[clip=true, trim={5cm, 0cm, 5cm, 0},  angle=90,  width=0.6\linewidth]{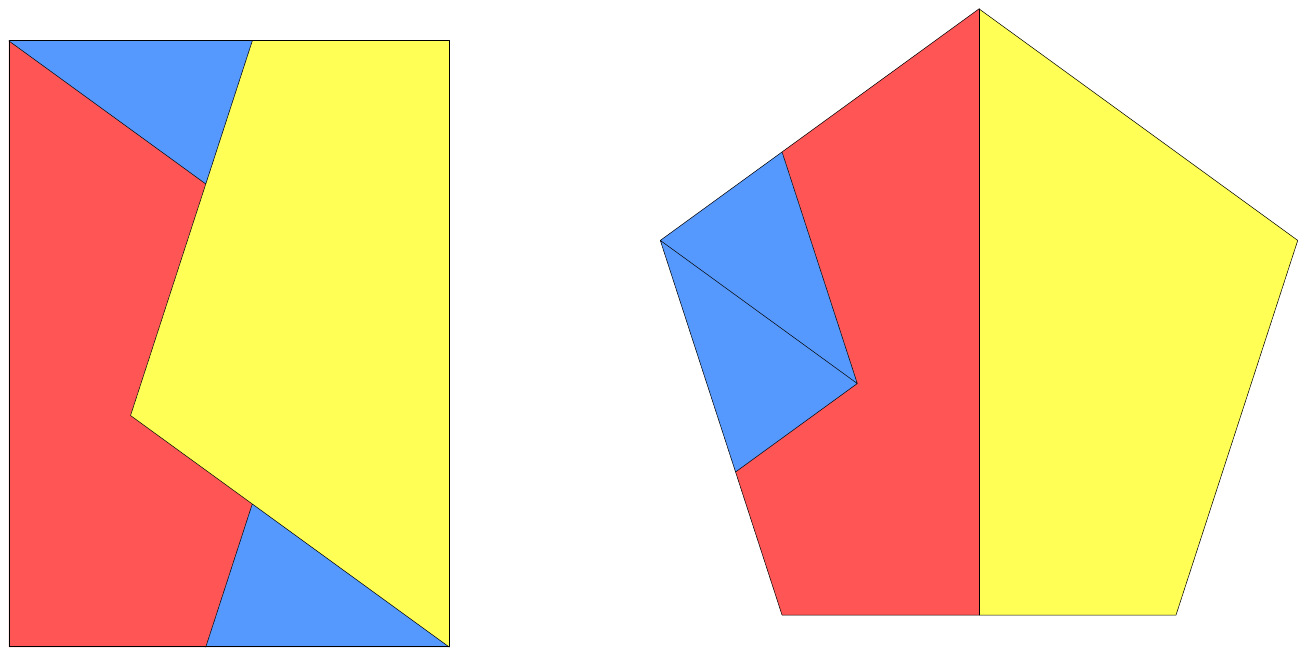}}
\caption{Gsellman's first pentagon dissection.}
\label{Fig5gona}
\end{figure} 


\subsection{Pentagon \#3}\label{Sec5C}
The first of Adam Gsellman's four-piece dissections of  a pentagon is shown
in Fig.~\ref{Fig5gona}. We do not know how Gsellman discovered it, but 
we have found that it can be obtained by a simple slide construction.
(It can also be obtained from a strip dissection.)

\begin{figure}[!ht]
\begin{center}
\begin{tikzpicture}[scale=2]
\def \eps1{0.02};
\def \SX{3};
\def \SY{.72};
\def \R{.851};
\def \w{1.118};
\def \h{1.539};
\coordinate (P1) at (0,\R);
\coordinate (P1a) at (0-\eps1,\R);
\coordinate (P1b) at (0+\eps1,\R);
\coordinate (P2) at (-.809,.263);
\coordinate (P3) at (-.5,-.688);
\coordinate (P4) at (.5,-.688);
\coordinate (P5) at (.809,.263);
\coordinate (Q1a) at (0-\eps1,-.688);
\coordinate (Q1b) at (0+\eps1,-.688);
\draw[ultra thick, brown] (P1a)--(P2);
\draw[ultra thick, brown] (P2)--(P3);
\draw[ultra thick, brown] (P3)--(Q1a);
\draw[ultra thick, brown] (Q1a)--(P1a);
\draw[ultra thick, cyan] (P1b)--(P5);
\draw[ultra thick, cyan] (P5)--(P4);
\draw[ultra thick, cyan] (P4)--(Q1b);
\draw[ultra thick, cyan] (Q1b)--(P1b);
\node at (-.35,0) {$\mathbf{A}$};
\node at (.4,0) {$\mathbf{B}$};
\coordinate (a1) at (\w+\SX, 0+\SY);
\coordinate (a2) at (0+\SX, 0+\SY);
\coordinate (a3) at (0+\SX, -\h+\SY);
\coordinate (a4) at (\w+\SX, -\h+\SY);
\coordinate (a5) at (.809+\SX, -.582+\SY);
\coordinate (b1) at (.618+\SX, 0+\SY);
\coordinate (b3) at (.309+\SX, -.951+\SY);
\coordinate (b5) at (0.5+\SX, -\h+\SY);
\draw[ultra thick, brown] (a2)--(a3);
\draw[ultra thick, brown] (a2)--(a5);
\draw[ultra thick, brown] (a3)--(b5);
\draw[ultra thick, brown] (b5)--(a5);
\draw[ultra thick, cyan] (b1)--(a1);
\draw[ultra thick, cyan] (a4)--(a1);
\draw[ultra thick, cyan] (b3)--(a4);
\draw[ultra thick, cyan] (b1)--(b3);
\node at (.2+\SX,-.45+\SY) {$\mathbf{A}$};
\node at (.9+\SX,-1.1+\SY) {$\mathbf{B}$};

\draw[thick] (a2)--(b1);
\draw[thick] (a5)--(b3);
\draw[thick] (b5)--(a4);
\node[above]  at (.3+\SX,0+\SY) {$\mathbf{s}$};
\node[below]  at (.55+\SX,-.77+\SY) {$\mathbf{s}$};
\node[below]  at (.82+\SX,-\h+\SY) {$\mathbf{s}$};
\end{tikzpicture}
\caption{A construction that produces the dissection in Fig.~\ref{Fig5gona} (see text for details).}
\label{Fig5gonc}
 \end{center}
\end{figure}
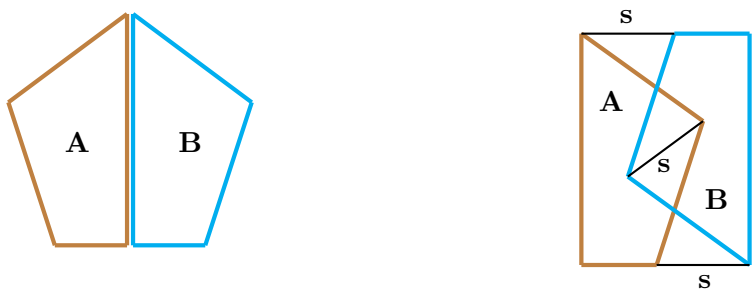

Cut the pentagon down the middle into two quadrilaterals $A$ and $B$, reflect $A$ in a vertical mirror, and rotate $B$ by $180^\circ$ (see Fig.~\ref{Fig5gonc}).
Now slide the pieces towards each other until they overlap in a parallelogram 
whose diagonal is equal in length to the gap in the top (and the bottom) edge.
Cut the parallelogram in the $A$ piece into  two equal isosceles triangles which can be rotated to complete the rectangle.

\ifnum \value{VarX}<1
{
\begin{figure}[!ht]
\centerline{\includegraphics[angle=90, width=6in]{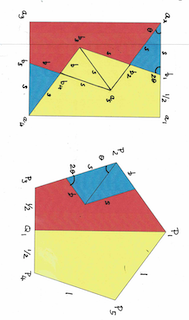}}
\caption{Names for points and lengths in dissection of Fig.~\ref{Fig5gona}.}
\label{Fig5gond}
\end{figure}
}
\else
{
\begin{figure}[!ht]
\centerline{\includegraphics[angle=90, width=6in]{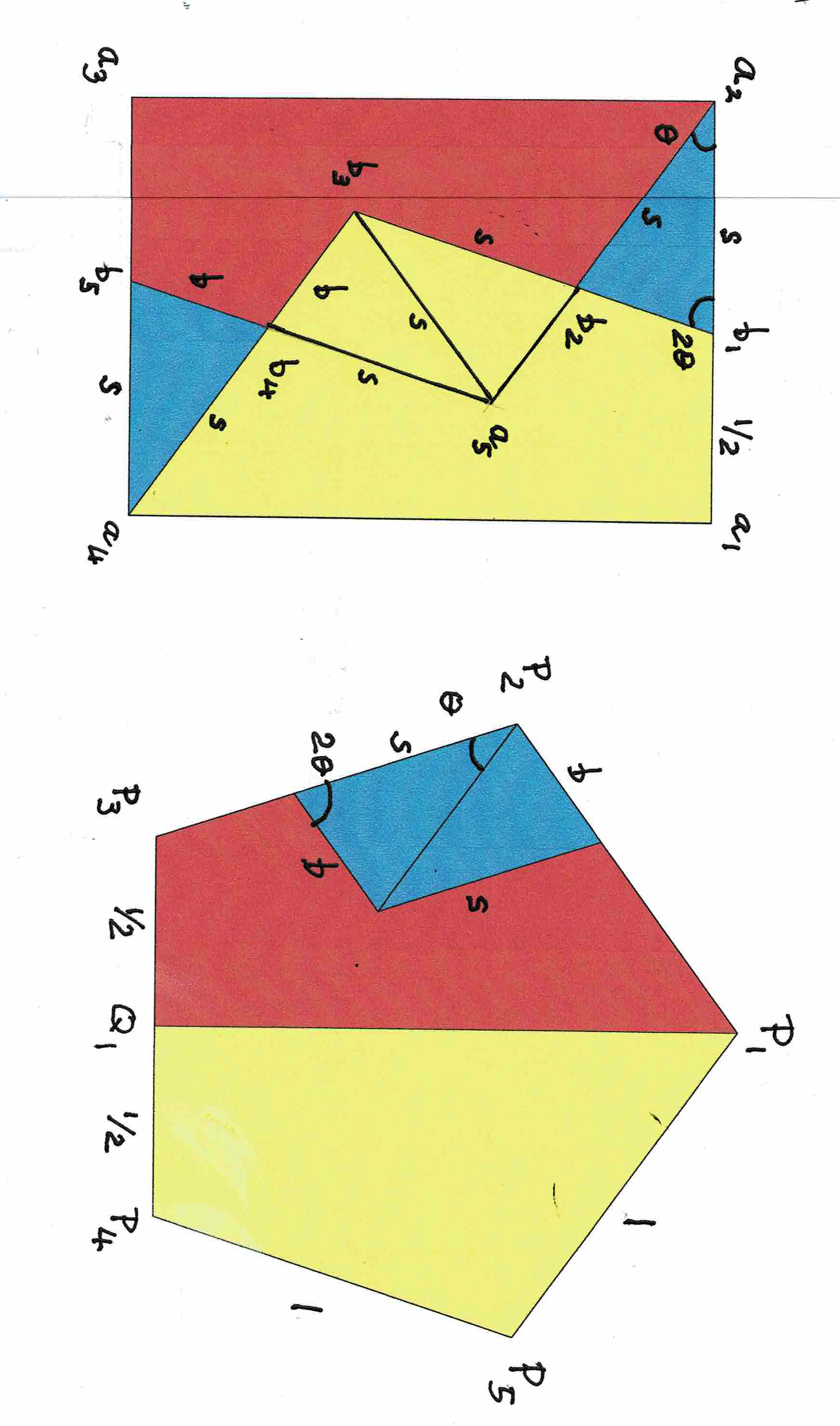}}
\caption{Names for points and lengths in dissection of Fig.~\ref{Fig5gona}.}
\label{Fig5gond}
\end{figure}
}
\fi

For the proof that this dissection is correct we label
the points as in Fig.~\ref{Fig5gond}.
$Q_1$ is the midpoint of the side $P_3 P_4$.
The key parameters are $\theta = \pi/5$, $R = 1/(2 \sin \theta)$,
and $d = R \cos \theta  = (\sqrt{5} + 1)^{3/2} /(4 \sqrt{2}\,5^{1/4})$.
The pentagon has height $h = R+d = (5^{1/4}/(4 \sqrt{2})) (\sqrt{5}+1)^{3/2}$.
This is also the height of the final rectangle, which (since we know the area)
has width $w = \sqrt{5}/2$.

We label the vertices of the rectangle 
 $a_1, \ldots, a_4$, and let $b_1, \ldots, b_5$, $a_5$  denote
 the points in the center of the figure. Finally, let $s$ and $b$ denote the side and base of the small isosceles triangles.
 
 Note that $a_1 b_1 b_3 a_4$ is a rotated copy of 
 $Q_1P_4 P_5 P_1$, and $a_2 a_3 b_5 a_5$ is a reflected copy of $P_1 P_2 P_3 Q_1$.
 In particular, $|a_1 b_1| = 1/2$, so $s = w - 1/2 = (\sqrt{5}-1)/2$.
 The base of the isosceles triangle is therefore $b = (3 - \sqrt{5})/2$.
 This implies $s+b = 1$, and we can now check that all the pieces fit together correctly.


\subsection{Pentagon \#4}\label{Sec5D}
\begin{figure}[!ht]
\centerline{\includegraphics[clip=true, trim={4cm, 0cm, 4cm, 0},  angle=90,  width=0.6\linewidth]{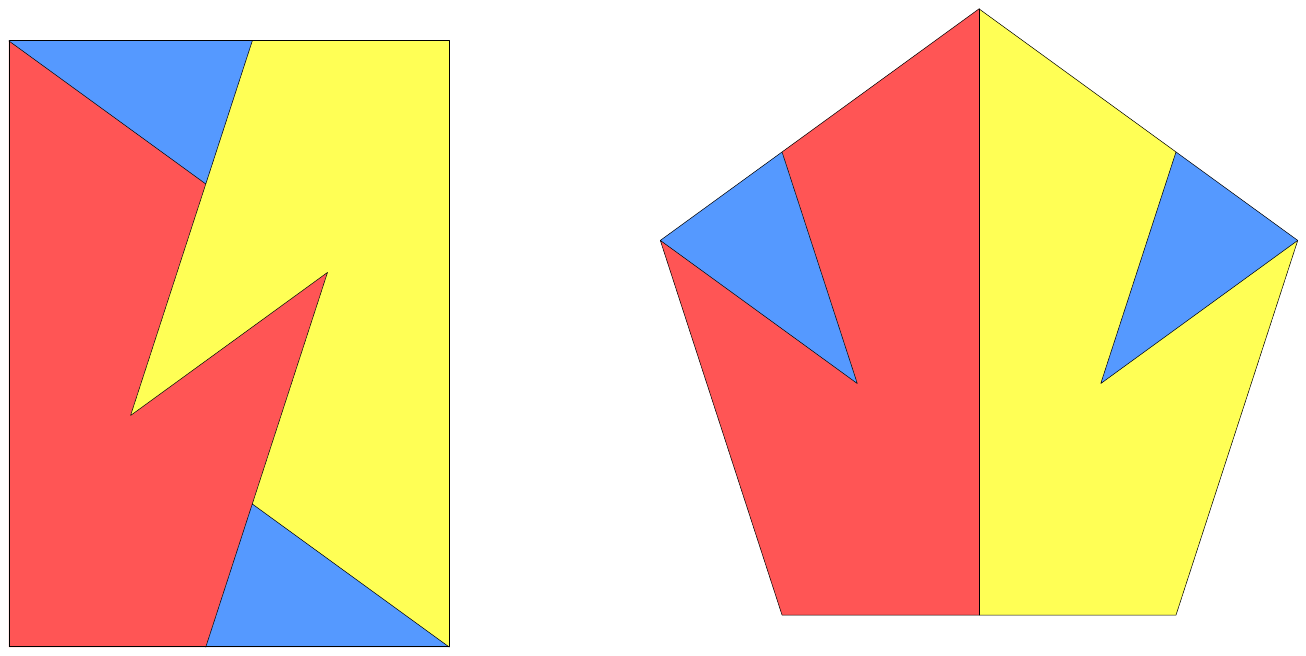}}
\caption{Gsellman's second pentagon dissection. }
\label{Fig5gonb}
\end{figure}

Gsellman's second pentagon dissection is shown in Fig.~\ref{Fig5gonb}.
It can be obtained by a similar slide construction.

\begin{figure}[!ht]
\centerline{\includegraphics[clip=true, trim={5cm, 0cm, 5cm, 0},  angle=90,  width=0.6\linewidth]{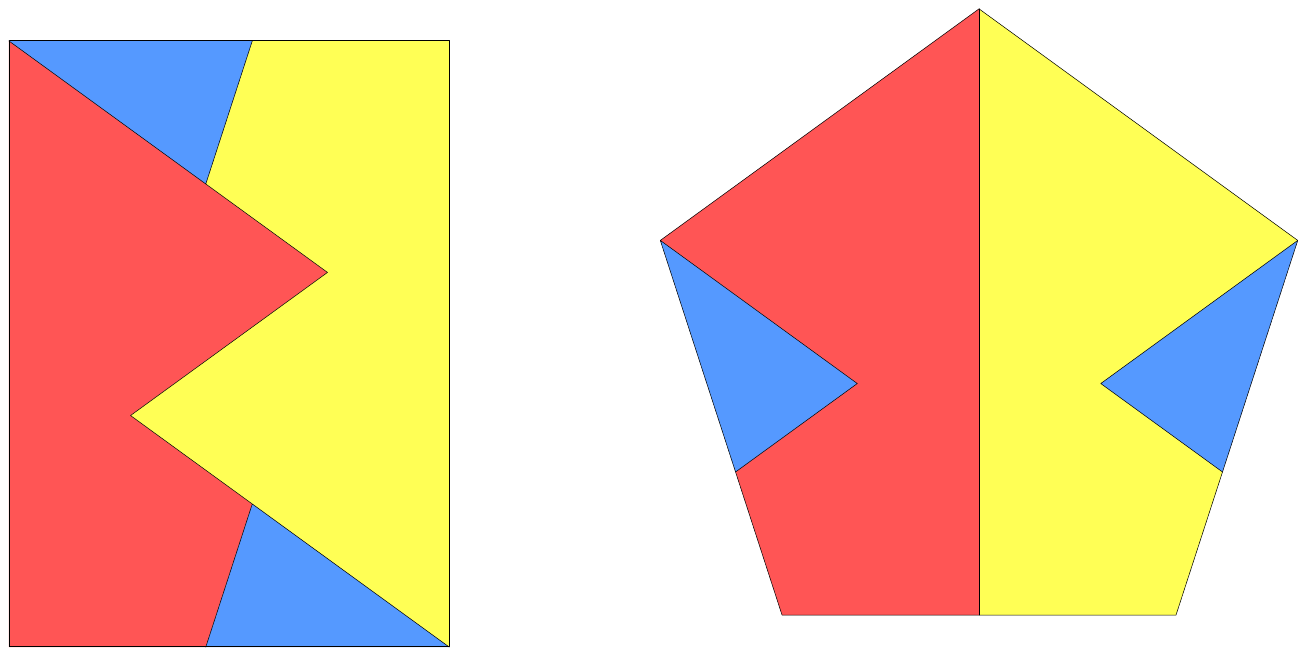}}
\caption{Another version of Fig.~\ref{Fig5gonb}. }
\label{Fig5gonbb}
\end{figure}

There is a third version of Gsellman's dissection which has the zigzag cut through the
center of the rectangle going the other way (Fig.~\ref{Fig5gonbb}). These are elegant dissections, 
but all three have the slight defect of requiring a piece to be turned over.        
 
\section{A five-piece dissection of a heptagon}\label{Sec7}

Even the great master Harry Lindgren \cite{Lind64} could only show that $s(7) \le 9$, but around 1995 
G.A.T.\  found a seven-piece dissection of a heptagon into a square, and it is  reasonable to conjecture that 
indeed $s(7) = 7$.
This dissection is described in \cite[pp.\ 128--129]{Fred97} and \cite{GDDb}.

\begin{figure}[!ht]
\centerline{\includegraphics[clip=true, trim={6cm, 0cm, 6cm, 0},  angle=90,  width=0.6\linewidth]
{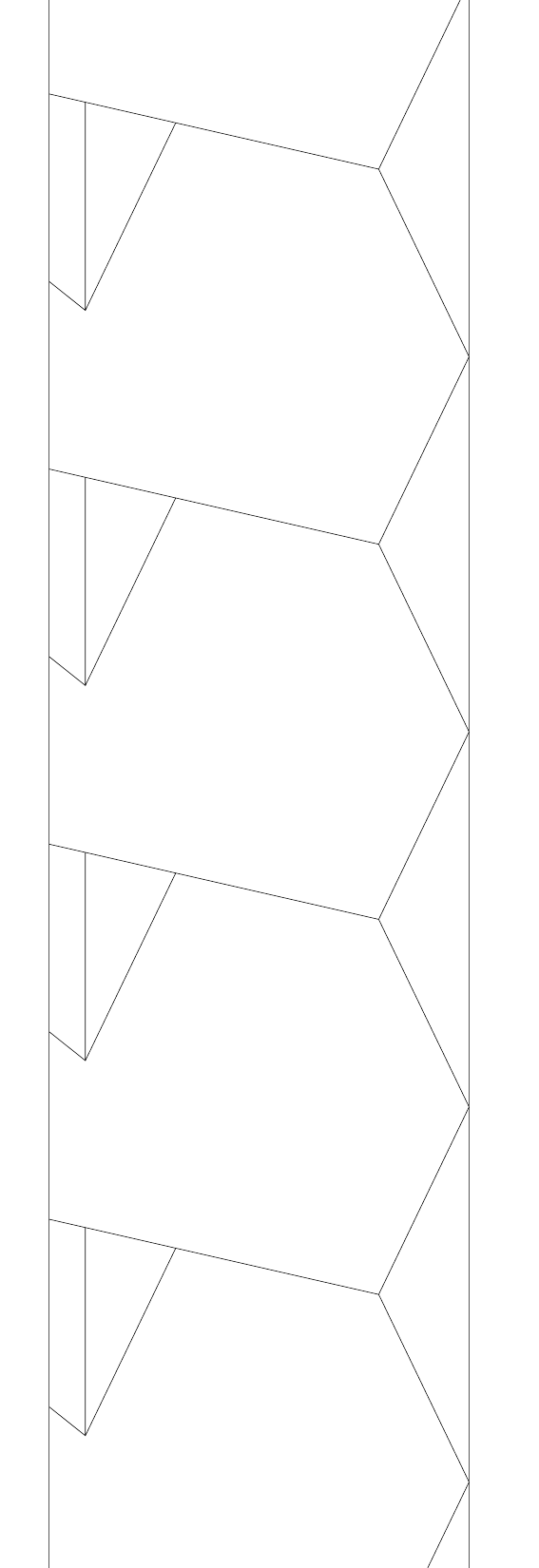}}
\caption{A heptagon strip used for the heptagon to rectangle dissection.}
\label{Fig7gonStripB}
\end{figure} 

\begin{figure}[!ht]
\centerline{\includegraphics[angle=0, width=4.5in]{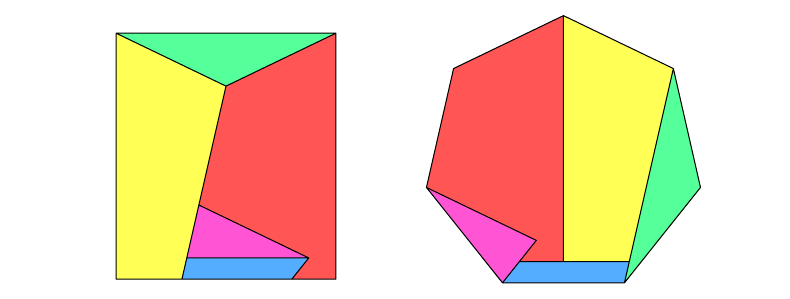}}
\caption{A five-piece dissection of a heptagon to a rectangle.}
\label{Fig7gon.d}
\end{figure}

For the rectangle problem, we start from a heptagon strip
(shown in Fig.~\ref{Fig7gonStripB})  that is a modification of
a strip used in the square case (compare~\cite[Fig.\ 11.28]{Fred97}).
By cutting a rectangle from this strip, we obtain a five-piece dissection of a heptagon to a rectangle, shown in Fig.~\ref{Fig7gon.d}.

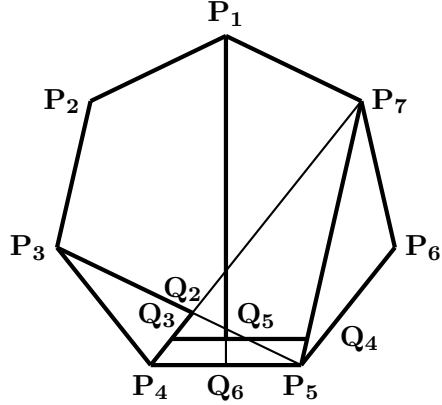
\begin{figure}[!ht]
\begin{center}
\begin{tikzpicture}[scale=2.0]
\coordinate(P1) at (0,1.152);
\coordinate(P7) at (.900, .718);
\coordinate(P2) at (-.900, .718);
\coordinate(P6) at (1.123, -.256);
\coordinate(P3) at (-1.123, -.256);
\coordinate(P5) at (.5, -1.038);
\coordinate(P4) at (-.5, -1.038);
\coordinate(Q2) at (-.2225, -.690);
\coordinate(Q3) at (-.3612, -.864);
\coordinate(Q6) at (.0, -1.038);
\coordinate(Q5) at (.0, -.8642);
\coordinate(Q4) at (.5397, -.8642);
\draw[ultra thick] (P1) -- (P2);
\draw[ultra thick] (P2) -- (P3);
\draw[ultra thick] (P3) -- (P4);
\draw[ultra thick] (P4) -- (P5);
\draw[ultra thick] (P5) -- (P6);
\draw[ultra thick] (P6) -- (P7);
\draw[ultra thick] (P7) -- (P1);

\draw[thick] (P7) -- (P4); 
\draw[thick] (P1) -- (Q6);
\draw[thick] (P3) -- (P5);
\draw[ultra thick] (P5) -- (P7);
\draw[ultra thick] (P3) -- (Q2);
\draw[ultra thick] (P4) -- (Q2);
\draw[ultra thick] (Q3) -- (Q4);
\draw[ultra thick] (P1) -- (Q5);

\node[above] at (P1) {$\mathbf{P_1}$};
\node[left] at (P2) {$\mathbf{P_2}$};
\node[left] at (P3) {$\mathbf{P_3}$};
\node[below] at (P4) {$\mathbf{P_4}$};
\node[below] at (P5) {$\mathbf{P_5}$};
\node[right] at (P6) {$\mathbf{P_6}$};
\node[right] at (P7) {$\mathbf{P_7}$};

\node[above] at (Q2) {$\mathbf{Q_2~~}$};
\node[above] at (Q3) {$\mathbf{Q_3}~~~$};
\node[right] at (Q4) {\hspace{.3cm}$\mathbf{Q_4}$};
\node[above] at (Q5) {\hspace{0.8cm}$\mathbf{Q_5}$};
\node[below] at (Q6) {$\mathbf{Q_6}$};

\end{tikzpicture}
\caption{Heptagon dissection (Fig.~\ref{Fig7gon.d}, right) showing labels for points.}
\label{Fig7gon.Pizza} 
 \end{center}
\end{figure}

Rather than following the path by which it was discovered, we will construct
this dissection directly from the heptagon, something that can
be done using only a straightedge and compass. 
The vertices of the $7$-gon are labeled $P_1, \ldots, P_7$ (see Fig.~\ref{Fig7gon.Pizza}).
Drop a perpendicular from $P_1$ to the midpoint $Q_6$ of $P_4  P_5$,
and draw chords $P_3 P_5$, $P_4 P_7$, and $P_5 P_7$.
Let $Q_2$ be the intersection of $P_3 P_5$ and $P_4 P_7$,
and let $Q_3$ be the midpoint of $P_4 Q_2$.
Draw a line $Q_3 Q_5 Q_4$ through $Q_3$ parallel to $P_4 P_5$.

Using these these lines as guides, we get the five pieces
by cutting $P_7 P_5 P_6 P_7$;
$P_1 Q_5 Q_4 P_7 P_1$;
$P_1 P_2 P_3 Q_2 Q_3 Q_5 P_1$;
$P_3 P_4 Q_2 P_3$;
and $Q_3 P_4 P_5 Q_4 Q_3$.

These pieces can then be rearranged to form the rectangle
as shown in Fig.~\ref{Fig7gon.d}, left.

To help anyone who wishes to verify the correctness of this dissection, we list some
key angles and lengths.
We set $\ta = \pi/7$ and note that $\cos (\ta) = .0.9009\ldots$ has minimal polynomial $8 x^3 - 4x^2 -4x +1$.
Then $\angle P_1 P_2 P_3 = 5 \theta$,
$\angle P_7 P_4 P_5 = 2\ta$,
$\angle P_4 P_3 P_5 = \angle P_5 P_7 P_6 = \ta$,
and $\angle P_4 P_5 P_7 = 4\ta$.
The chord $P_7 P_5$ has length $L_{2,7} = 2 \cos \ta$, and
$|P_4 Q_3| = |Q_3 Q_2| = 1/2 - (\sec \ta)/4$.
The trapezoidal piece has cross-section $|Q_5 Q_6| =  (4 \cos \ta -3)/(8 \sin \ta)$.
The rectangle has height $7 /(8 \sin \ta) $ and width $2 \cos \ta$.


\begin{figure}[!ht]
\centerline{\includegraphics[angle=0, width=4in]{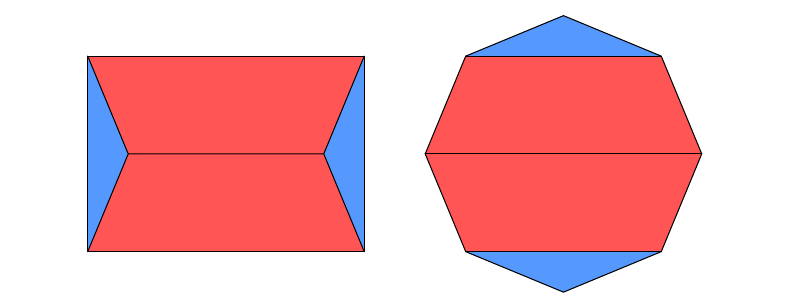}}
\caption{The classic four-piece  octagon to rectangle dissection \cite{GDDb}.}
\label{Fig8A}
\end{figure} 

\section{Four-piece dissections of an octagon}\label{Sec8}

It appears that five pieces are needed to
dissect an octagon into a square (see Fig.~\ref{FigOctSquare}),
whereas four pieces are enough if we only
want a rectangle (Fig.~\ref{Fig8A}).  The former has cyclic four-fold symmetry, 
while the latter has the symmetry of a Klein $4$-group.

\begin{figure}[!ht]
\centerline{\includegraphics[clip=true, trim={0cm, 0.5cm, 0cm, 0.5cm},  angle=0,  width=0.6\linewidth]{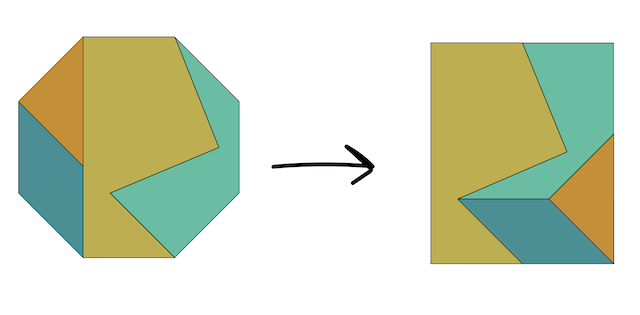}}
\caption{Gsellman's first  four-piece dissection of an octagon. See Fig.~\ref{Fig8C} for details.}
\label{Fig8B}
\end{figure}

\begin{figure}[!ht]
\centerline{\includegraphics[angle=0, width=5in]{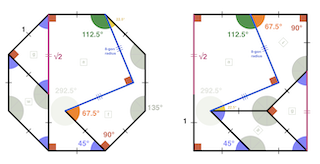}}
\caption{Gsellman's detailed description of the dissection in Fig.~\ref{Fig8B}.}
\label{Fig8C}
\end{figure}

\begin{figure}[!ht]
\centerline{\includegraphics[clip=true, trim={0cm, 0.5cm, 0cm, 0.5cm},  angle=0,  width=0.6\linewidth]{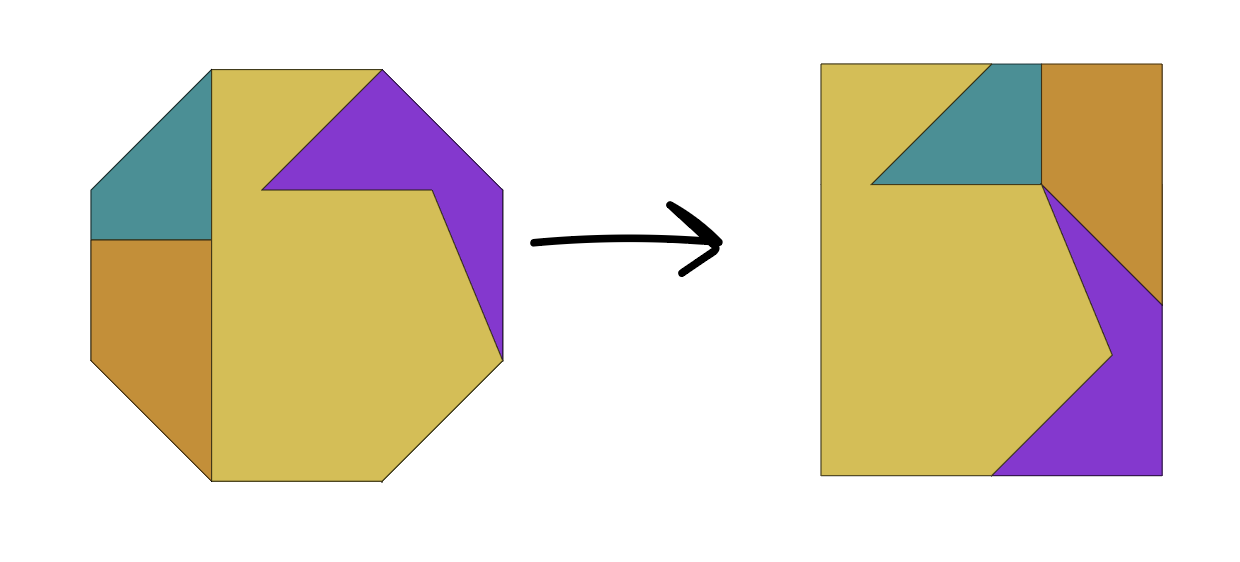}}
\caption{Gsellman's second four-piece dissection of an octagon.}
\label{Fig8D}
\end{figure}

Adam Gsellman \cite{Gsel23} has found two other four-piece dissections,
shown in Figs.~\ref{Fig8B}-\ref{Fig8D}.
The description in Fig.~\ref{Fig8C} is self-explanatory (the angles
are multiples of $\pi/8$ and the only irrationality needed is $\sqrt{2}$).
The second dissection (Fig.~\ref{Fig8D}) is very similar.

These two dissections are admittedly less elegant than that in Fig.~\ref{Fig8A},
and require pieces to be turned over,
but we include them because, as the example in Fig.~\ref{Fig10A} shows,
complicated non-convex dissections may be needed to get the minimum number of pieces.

\section{A  seven-piece dissection of a 9-gon}\label{Sec9}

\begin{figure}[!ht]
\centerline{\includegraphics[clip=true, trim={3cm, 0cm, 5cm, 0cm},  angle=90,  width=0.670\linewidth]{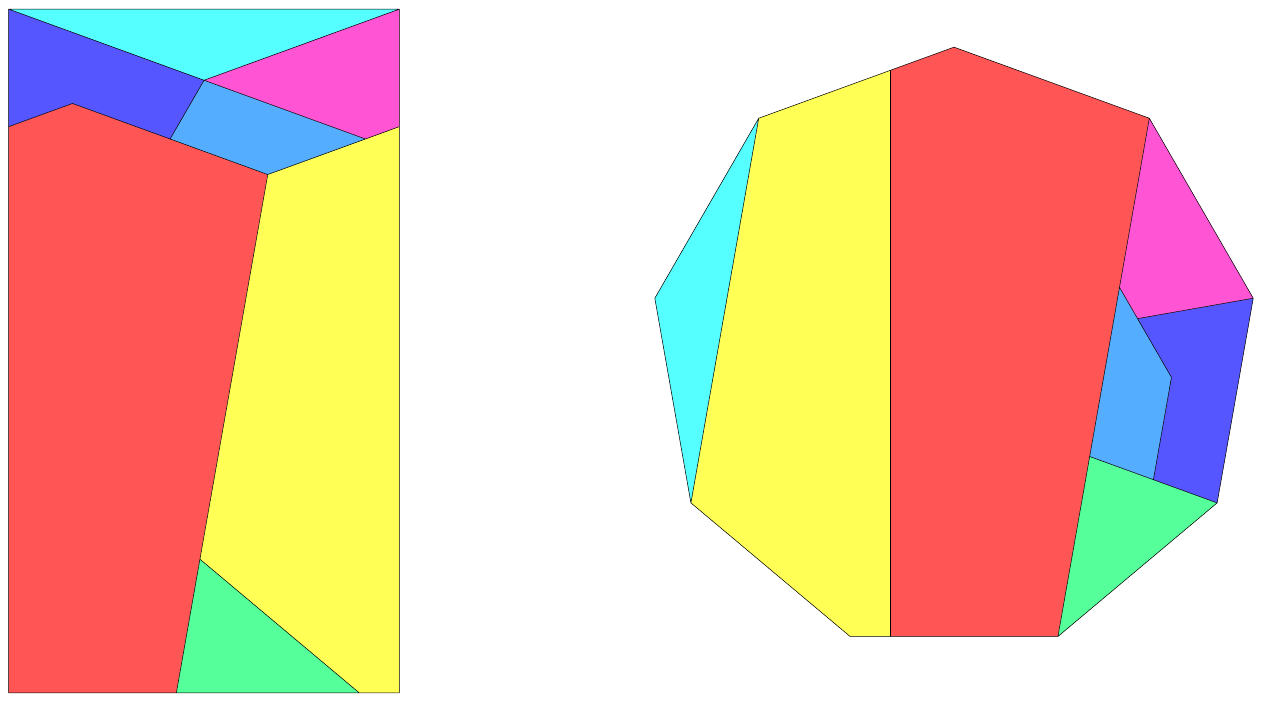}}
\caption{ A seven-piece dissection of a $9$-gon into a rectangle.}
\label{Fig9A}
\end{figure}

\begin{figure}[!ht]
\centerline{\includegraphics[clip=true, trim={3cm, 0cm, 6cm, 0},  angle=270,  width=0.6\linewidth]{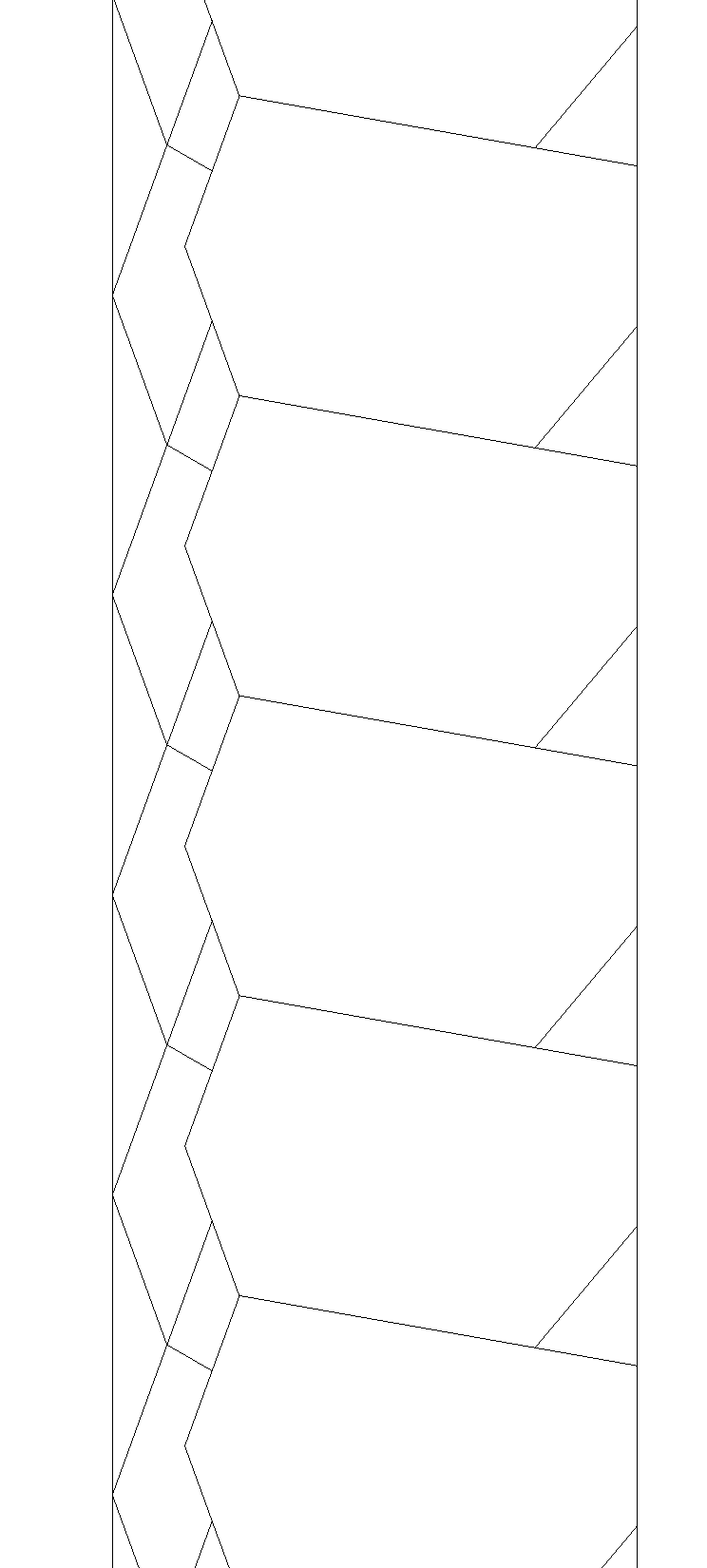}} 
\caption{ The strip dissection of a $9$-gon which led to Fig.~\ref{Fig9A}.}
\label{Fig9B}
\end{figure}

Figure~\ref{Fig9A} shows a seven-piece dissection of a $9$-gon into a rectangle,
which is two fewer pieces than the best dissection into a square presently known.
It was obtained from the strip dissection of a $9$-gon shown in Fig.~\ref{Fig9B},
by cutting a rectangle from the strip.

\begin{figure}[!ht]
\centerline{\includegraphics[clip=true, trim={0cm, 2cm, 0cm, 4cm},  angle=0,  width=0.6\linewidth] {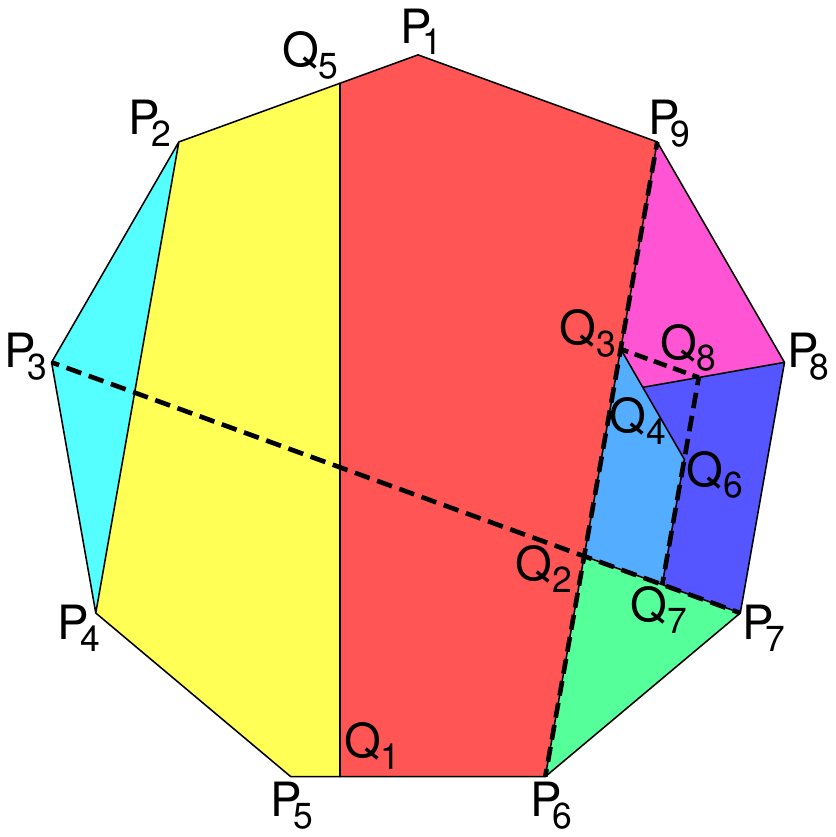}}
\caption{Labels for points in the $9$-gon.}
\label{Fig9C}
\end{figure} 

As with the heptagon,  we will construct 
this  dissection directly from the $9$-gon,  using only a straightedge and compass.

The vertices of the $9$-gon are labeled $P_1, \ldots, P_9$ (see Fig.~\ref{Fig9C}), and we
use the coordinates established in~\S\ref{Sec2}. 
Also $\theta = \pi/9$,  $C_1 = \cos \theta$ has minimal polynomial $8x^3-6x-1$,
and $S_1 = \sin \theta$ has minimal polynomial 
$64 x^6 - 96 x^4 + 36 x^2 - 3$.
(For this reason, we write our expressions as rational functions of $\cos \theta$,
with at most  linear terms in $\sin \theta$.)

To obtain the dissection we first draw chords $P_2-P_4$, $P_3-P_7$,
and $P_6-P_9$.
Then $Q_2$ is the intersection of $P_3-P_7$ and $P_6-P_9$, $Q_3$ is the midpoint of 
$Q_2-P_9$, and $Q_7$ is the midpoint of $Q_2-P_7$.
Draw a line segment $Q_7-Q_6$ of length $1/2$  parallel to $P_7-P_8$,
join $Q_3$ to $Q_6$, and locate $Q_4$ at the intersection of  $Q_3-Q_6$ 
and a perpendicular drawn from $P_8$ to the midpoint of $P_3-P_4$. Finally 
$Q_5$ is on $P_1-P_2$ at distance $|Q_4Q_6|$ from $P_1$, and 
$Q_5-Q_1$ is perpendicular to $P_5-P_6$.

To assist in the analysis we define a further point $Q_8$ 
at the intersection of 
$Q_7-Q_6$ (extended) and $Q_4-P_8$. Then $Q_4 Q_6 Q_8$ is an isosceles triangle
and $Q_3 Q_2 Q_7 Q_8$ is a parallelogram.

The seven pieces of the dissection can now be found
by making the cuts indicated by the colored regions on the right of
Fig.~\ref{Fig9A}.
To complete the proof that the dissection is correct, we must verify that
the pieces can be rearranged to form the rectangle on the left of Fig.~\ref{Fig9A}.
We will not take the space to do that here,
but to assist the reader we give two key lengths. 
The length 
\beql{Eq9.1}
|Q_2Q_7| = |Q_2 P_7|/2 =  |Q_4 Q_6| = |Q_4Q_8| = |Q_8Q_3| = |P_1 Q_5| = 
 |Q_2Q_3|-\frac{1}{2} =  \frac{C_1}{2 C_1+1} = 0.3263\ldots
\eeq
plays a central role, as does
$|P_8 Q_4| = 3/(8 S_1 (C_1+1)) = 0.5652$.
The rectangle has width $2 C_1 = 1.8793\ldots$ and height $9/(8 S_1) = 3.2892\ldots$.

Some of the equalities in \eqn{Eq9.1} are by no means obvious.
They do not follow directly from the geometry, but depend on the fact
that $\cos \theta$ satisfies a cubic equation. As discussed in Remark~\ref{RemSEC},
we can find exact expressions for the coordinates of the points.
We {\em could} do this by solving the appropriate equations,
using a computer algebra system such as Maple,
but we have found it a lot easier to use another computer algebra system, WolframAlpha,
and ask it to find the coordinates for us.

For example, the first step in finding the present dissection is to find $Q_2$.
Using Maple, and working to $20$ decimal places, we find that
$$
Q_2 ~=~ (0.65270364466613930216, -0.50771330594287249271)\,.
$$
We now ask WolframAlpha  to express these two numbers in terms of 
$\cos \theta$, $1/\cos \theta$, $\sin \theta$, and $1/\sin \theta$.
The result (setting $C_1 = \cos \theta$, $S_1 = \sin \theta$) is
$$
Q_2 ~=~ \left(\frac{2 C_1}{2 C_1 +1},~ -(2 S_1 - \frac{1}{S_1} +\sqrt{3})\right)\,.
$$
We do this for all the points. Another example is
\begin{align} 
Q_8 & ~=~ (1.1028685319524432095, 0.19446547835755153996) \nonumber \\
& ~=~ \left( \frac{C_1}{2} + \frac{1}{8 C_1} + \frac{1}{2}.
~\frac{1}{8 S_1} - \frac{S_1}{2} \right)\,.
\end{align}

The verification of the equalities in \eqn{Eq9.1}
is then a routine calculation (using  Maple's \texttt{simplify(..., trig);}  command).


\section{A  four-piece dissection of a 10-gon}\label{Sec10}

In the final appendix (``Recent Progress'')  to his 1964 book \cite{Lind64},
Lindgren gave a new strip based on the $10$-gon (shown in Fig.~\ref{Fig10Bpre} below),
and used it to obtain several new dissections, including an eight-piece
dissection of a $10$-gon to a square, and a seven-piece dissection to a golden rectangle.
As Frederickson reports in  \cite[Ch.~11]{Fred97}, G.A.T.\ was then able to show 
that in the dissection to a square, 
two of Lindgren's pieces could be merged,
leading to a seven-piece dissection to a square,  still the record.
This dissection is also described in the {\em Variable Strips} section of~\cite{GDDb}.
If we draw vertical lines across Lindgren's strip (without changing it),
we obtain a five-piece dissection of a $10$-gon
to a (non-golden) rectangle, as shown in Fig.~\ref{Fig10Lind}.
There is a small range of possibilities for the positions of these vertical lines.
In Fig.~\ref{Fig10Lind} they have been placed in the middle of their range, in order to 
obtain the most symmetric dissection.

\begin{figure}[!ht]
\centerline{\includegraphics[angle=0, width=4.5in]{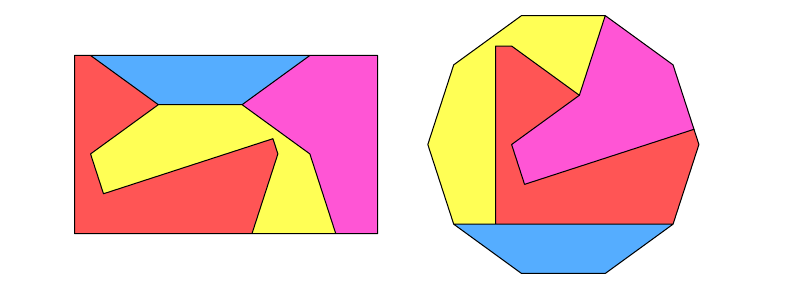}}
\caption{ A four-piece dissection of a 10-gon into a rectangle.}
\label{Fig10A}
\end{figure}

Remarkably, if the goal is only to obtain a rectangle, it is possible to modify
Lindgren's strip  (Fig.~\ref{Fig10Lind}) to
get a four-piece dissection.  The modified strip is shown in Fig.~\ref{Fig10B},
and the dissection itself  in Fig.~\ref{Fig10A}.
To go from Fig.~\ref{Fig10Lind} to Fig.~\ref{Fig10B} we merge the two right-most 
pieces of the rectangle, forming a church-shaped piece, and compensate by 
dividing the large piece into two by a zig-zag cut.

\begin{figure}[!ht]
\centerline{\includegraphics[clip=true, trim={7cm, 0cm, 7cm, 0},  angle=90,  width=0.6\linewidth]
{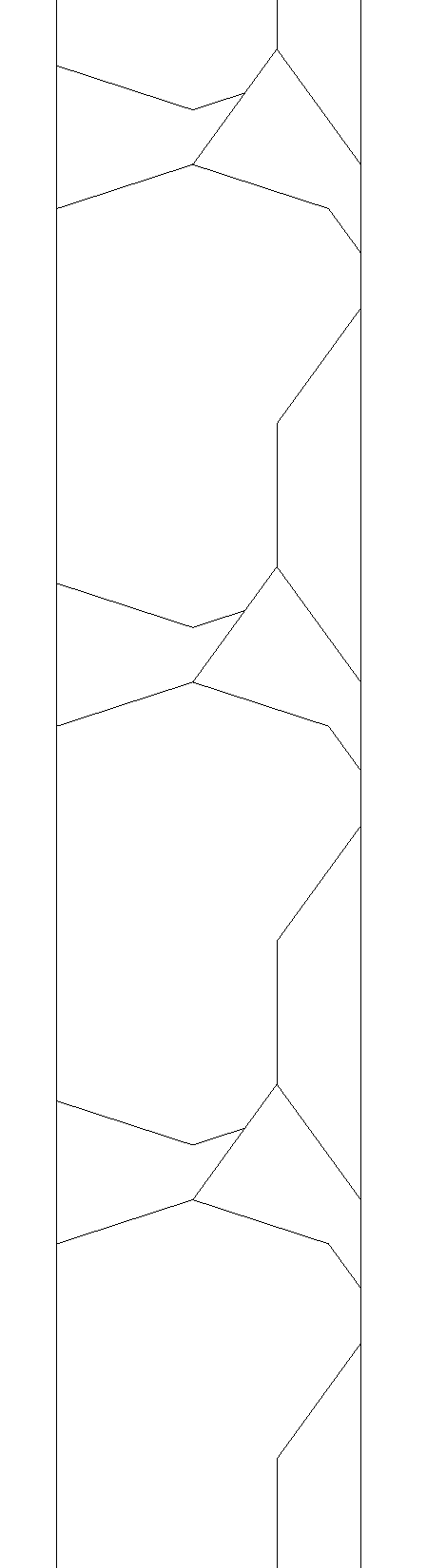}}
\caption{Lindgren's 1964 $10$-gon strip \cite{Lind64}.}
\label{Fig10Bpre}
\end{figure} 

\begin{figure}[!ht]
\centerline{\includegraphics[clip=true, trim={7cm, 0cm, 7cm, 0},  angle=90,  width=0.6\linewidth]{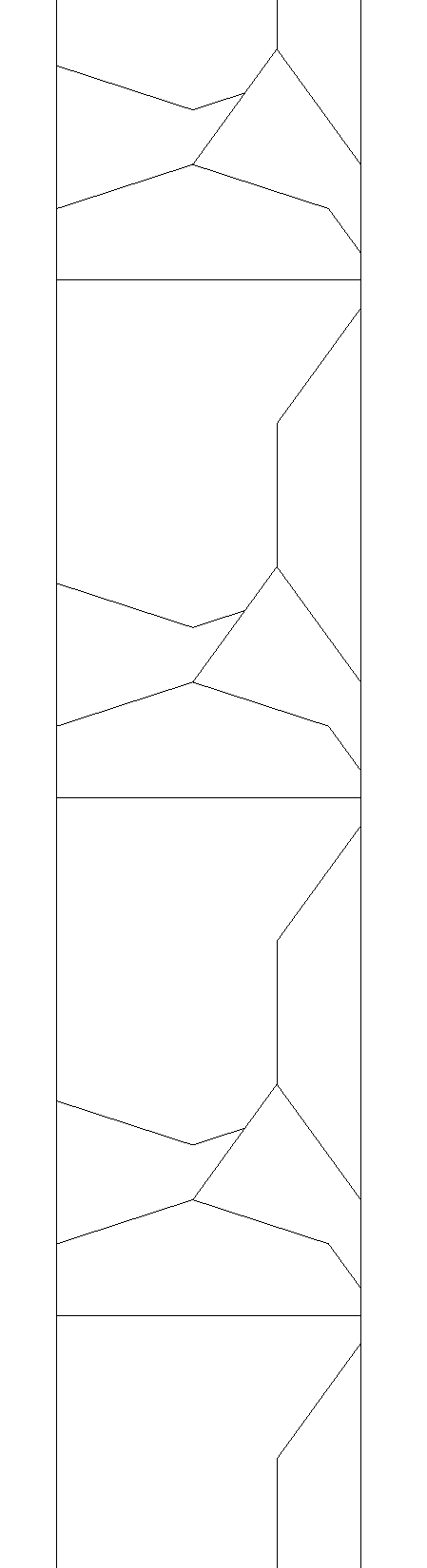}}
\caption{A five-piece $10$-gon to rectangle dissection obtained from Fig.~\ref{Fig10Bpre}.}
\label{Fig10Lind}
\end{figure}

\begin{figure}[!ht]
\centerline{\includegraphics[clip=true, trim={7cm, 0cm, 7cm, 0},  angle=90,  width=0.6\linewidth]{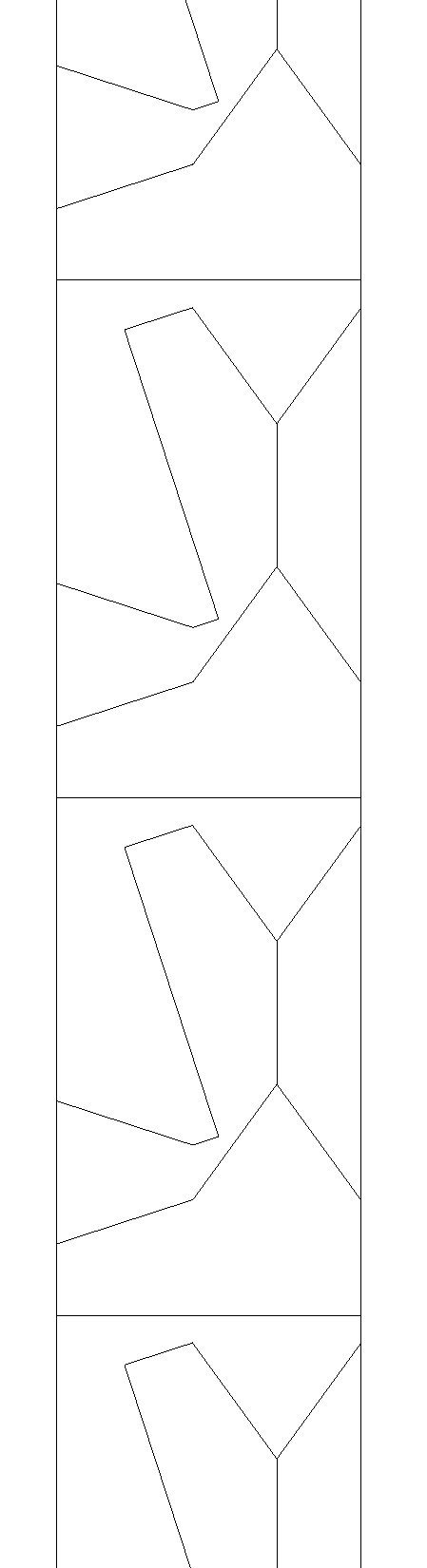}}
\caption{The strip which gives the four-piece dissection of the $10$-gon to a rectangle.}
\label{Fig10B}
\end{figure} 

This is one of the most complicated dissections in the article, and we give
a precise straightedge and compass construction starting from  
the rectangle.in Fig.~\ref{Fig10B}.
 
We first construct an intermediate rectangle with five pieces, and then shift
it slightly to save a piece. The angles involved are $\theta = \pi/10 = 18^{\circ}$, $2 \theta$,
and $\phi = 4 \theta$.

\begin{figure}[!ht]
\centerline{\includegraphics[angle=0, width=4.3in]{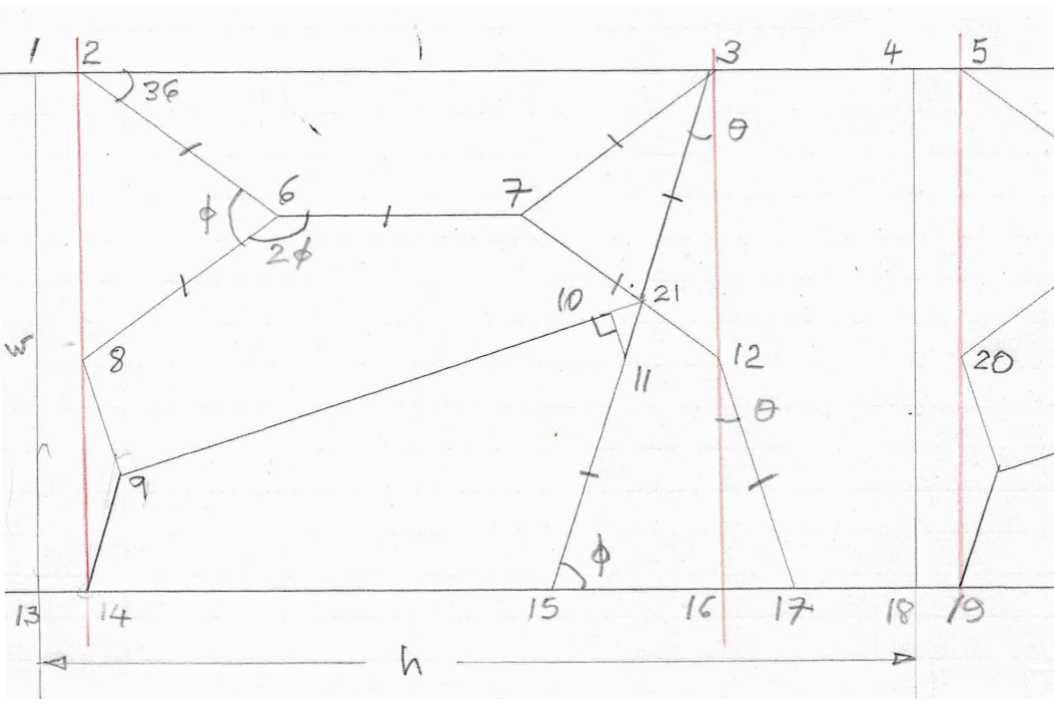}}
\caption{Labels for points used to construct the four-piece dissection of the $10$-gon.}
\label{Fig10C}
\end{figure}

The intermediate rectangle has vertices labeled $2, 14, 19, 5$ in Fig.~\ref{Fig10C};
the final rectangle has vertices $1, 13, 18, 4$. We place
the origin of coordinates for the rectangle  near the bottom left corner,  at the point $14 = (0,0)$. 
The $10$-gon has area $\frac{5}{2 \tan \theta}$ (see \eqn{EqArea}), 
and we take the width of the strip to be $w = \sqrt{5}\, \cos \theta$.
The other dimension of the rectangle is its height $h = 2 \sqrt{5}\,\cos 2\theta$.
(After a series of relabelings, the rectangle as drawn in Fig.~\ref{Fig10B} has ended up
with height $w$ and width $h$. We hope the reader will forgive us!)

The coordinates of the points $2, 14, 17, 5$ are therefore $(0,w)$, $(0,0)$, $(h,0)$,
and $(h.w)$.
We draw a network of lines as follows.  Starting at point $2$, we draw line
segments of length $1$ from $2$ to $6$ to $7$ to $3$,
and from $6$ to $8$, $3$ to $21$,  and $7$ to $12$ to $17$ to $15$ to $11$  (the angles are indicated in the figure).
We then complete the line $3$ to $15$. We also draw line segments of length $1/2$ from $8$ to $9$ to $14$.
For the two final lines we join $9$ to $21$ and draw the perpendicular from $11$ to $10$.
The coordinates have been chosen so that  several coincidences occur.\footnote{Similar to those in \eqn{Eq9.1}, but less dramatic, since $\sin \pi/10$ is only a quadratic irrationality.}
The points $8$, $11$, $12$, and $20$ (in the adjacent rectangle in the strip) are collinear. Also 
the distance from $9$ to $10$ turns out to be equal to $w$. The angle $\angle 8, 9, 21$ is a right angle.
The central point $21$ has coordinates 
$(2 + \sin \theta, 2 \sin 2\theta)$.
The distance from $10$ to $11$ is $x := \frac{3-\sqrt{5}}{4}$,
and we get the final rectangle by shifting the intermediate rectangle to the left by that amount.

We get the four pieces in the dissection as follows. 
The quadrilateral piece (the ``dish'') is obtained by cutting along the path $2, 6, 7, 3, 2$.
For the hexagon (the ``church''), cut along $3, 7, 12, 17, 18, 4, 3$.
For the first $9$-gon (the ``hammer"), cut along $6, 8, 9, 10, 11, 15, 17, 12, 7, 6$,
and for the second $9$-gon (the ``triangle''), cut along  $2, 1, 13, 15, 11, 10, 9, 8, 6, 2$.
By moving the edge of the rectangle to the left so that it no longer passes through the point $8$ we
have  reduced the number of pieces from five to four.


\begin{figure}[!ht]
\centerline{\includegraphics[angle=0, width=5in]{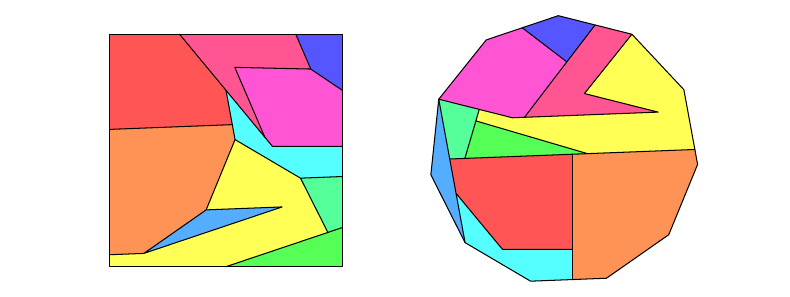}}
\caption{A ten-piece dissection of an $11$-gon into a square \cite{GDDb}.}
\label{Fig11A}
\end{figure} 

\section{Dissecting  an 11-gon to a square and to a rectangle}\label{Sec11}

\subsection{A ten-piece dissection of an 11-gon to a square}\label{Sec11s}

\begin{figure}[!ht]
\centerline{\includegraphics[angle=90, width=4.3in]{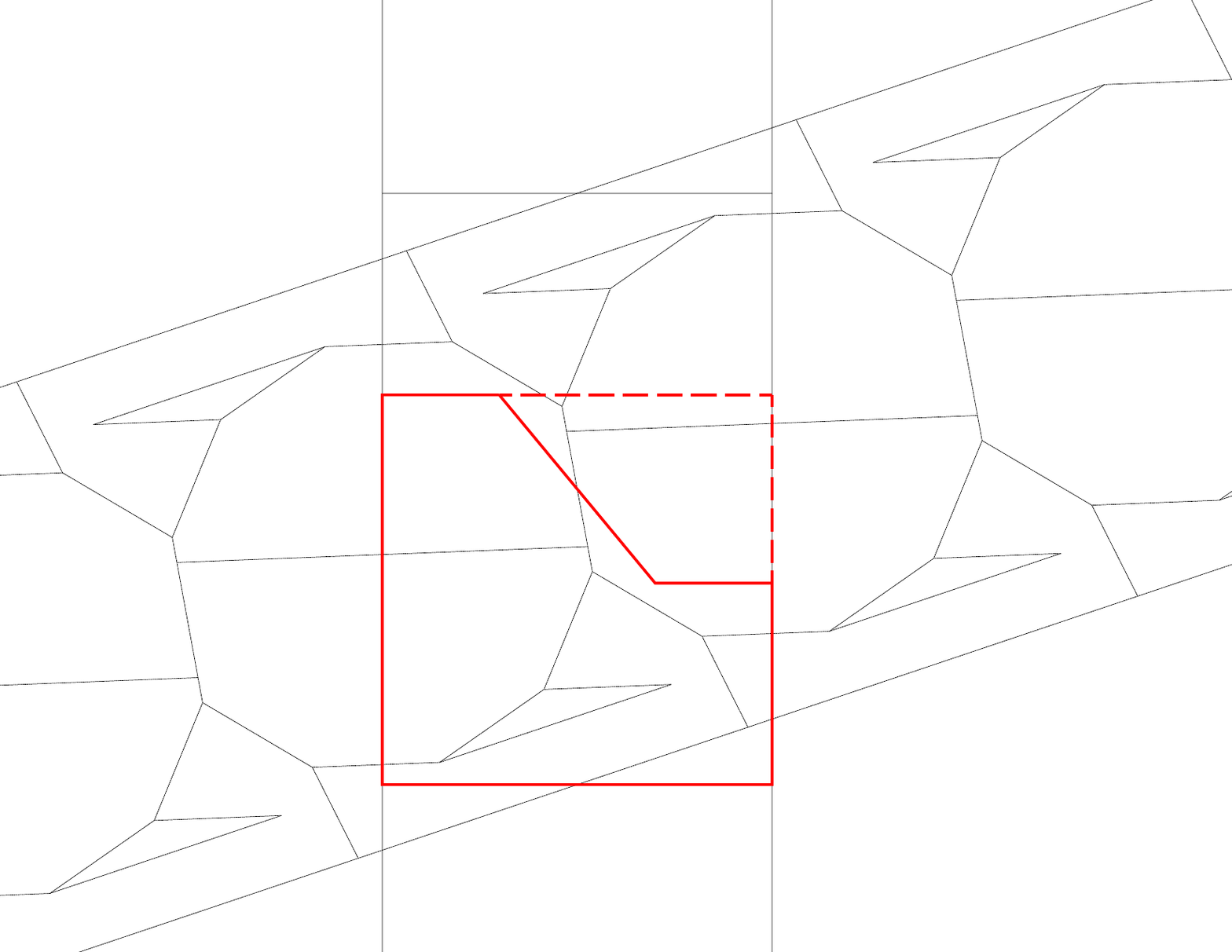}}
\caption{First of two strip superpositions used in construction of the $11$-gon to square dissection of Fig.~\ref{Fig11A}.}
\label{Fig11B}
\end{figure}

\begin{figure}[!ht]
\centerline{\includegraphics[clip=true, trim={5cm, 0cm, 5cm, 0cm},  angle=90,  width=0.8\linewidth]{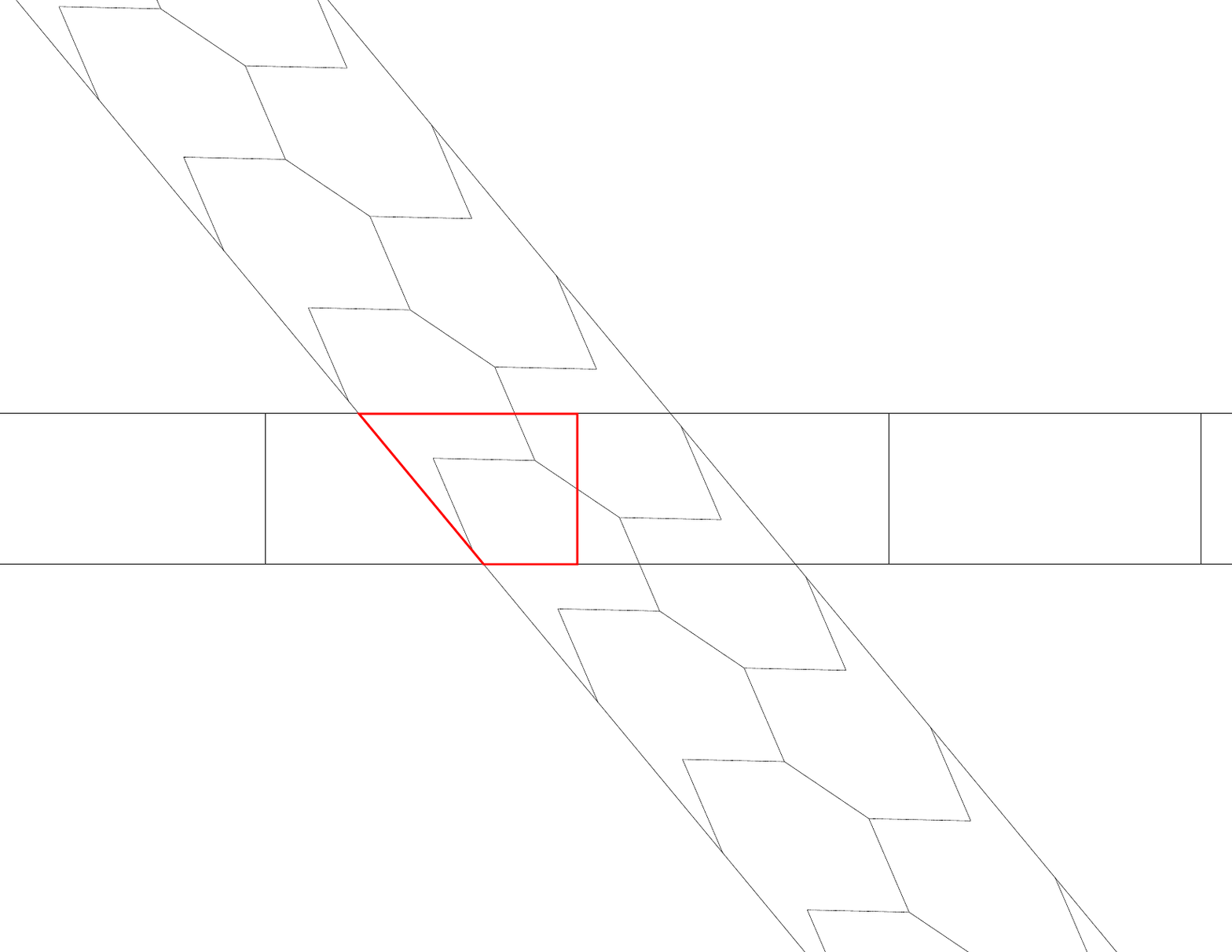}}
\caption{Second of two strip superpositions used in construction of the $11$-gon to square dissection of Fig.~\ref{Fig11A}. The large convex piece in the red region  has one very short edge, and is actually a hexagon.
This hexagon (colored pink) appears in both illustrations in Fig.~\ref{Fig11A}. }
\label{Fig11C}
\end{figure}

Before the appearance of \cite{GDDb} there had been little work on dissections of the $11$-gon:
this polygon is not mentioned in any of \cite{Fred97, Lind64, Lind72}.
G.A.T.'s  ten-piece dissection of an $11$-gon into a square was given in \cite{GDDb},
and was described by Frederickson in  \cite{Fred97a, Fred06}.
It is shown here in  Fig.~\ref{Fig11A}. 
It can be  obtained by taking the $11$-gon and constructing  the two superpositions of strips shown
in Figs. \ref{Fig11B} and \ref{Fig11C}. 
When the quadrilateral outlined in red in Fig.~\ref{Fig11C}
is combined with the hexagon outlined in red in Fig.~\ref{Fig11B},
the result is the dissected square on the left of Fig.~\ref{Fig11A}.

\subsection{Our first  nine-piece dissection of an 11-gon to a rectangle}\label{Sec11r}

A piece can be saved if our goal is only to dissect the $11$-gon into a rectangle.
We have found several examples of nine-piece $11$-gon to rectangle
dissections, two of which are described here and in the next  section.
Our first construction is similar to the $11$-gon to square dissection of \S\ref{Sec11s}.
The proof of correctness involves an interesting interplay between the two superpositions.
A second reason for including this proof is that similar arguments can be used to give
a proof of the $11$-gon to square dissection mentioned above.

We start from the $11$-gon and construct  two superpositions of strips 
(see Figs.\ \ref{Fig11D} and \ref{Fig11Db}), and 
when the quadrilateral outlined in red in Fig.~\ref{Fig11D}
is combined with the hexagon outlined in red in Fig.~\ref{Fig11Db},
the result is the dissected rectangle on the left of Fig.~\ref{Fig11E}.

\begin{figure}[!ht]
\centerline{\includegraphics[angle=6, width=2.5in]{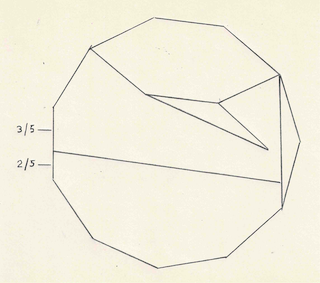}}
\caption{A 5-piece dissection of an $11$-gon used to form the superposition of Fig.~\ref{Fig11D}.}
\label{Fig11S1}
\end{figure}

To build these superpositions, we first cut the $11$-gon into five pieces, as shown in 
Fig.~\ref{Fig11S1}. The sides of the $11$-gon have length $1$, as usual,  and the two long skinny triangles have
a base of length $L_{11,2} = 2 \cos \theta$, where $\theta = \pi/11$.
The long cut through the middle of Fig.~\ref{Fig11S1} begins at a point 
three-fifths of the way along an edge. This parameter could be changed, but $3/5$ seems to be a good choice.
This cut is parallel to the top edge of the polygon, and has length 
$L_{11,4} = 2(\cos \theta + \cos 3\theta)$.

The heptagonal tadpole-shaped structure at the top of Fig.~\ref{Fig11S1} was the repeating element in 
Fig.~\ref{Fig11C} and will be used again in the
second superposition (Fig.~\ref{Fig11Db}).
It can be formed from three pieces cut from the $11$-gon along chords, as follows:
\begin{center}
\begin{tikzpicture}[scale=1.0]
\coordinate(C) at (0.5, 1.703);
\coordinate(B) at (-0.5, 1.703);
\coordinate(F) at (-0.5, .621);
\coordinate(E) at ( 0.5, .621);
\coordinate(D) at (1.341, 1.162);
\coordinate(G) at (-1.341, 2.243);
\coordinate(A) at (-1.341, 1.162);
\draw[ultra thick] (A) -- (B); 
\draw[thick] (B) -- (C);
\draw[ultra thick] (C) -- (D);
\draw[ultra thick] (A) -- (F);
\draw[ultra thick] (F) -- (E);
\draw[ultra thick] (E) -- (D);
\draw[ultra thick] (G) -- (B);
\draw[ultra thick] (G) -- (C);
\draw[thick] (A) -- (D);
\node[left] at (A) {$\mathbf{A}$};
\node[left] at (B) {$\mathbf{B}~$};
\node[right] at (C) {$~\mathbf{C}$};
\node[right] at (D) {$\mathbf{D}$};
\node[right] at (E) {$~\mathbf{E}$};
\node[left] at (F) {$\mathbf{F}~$};
\node[left] at (G) {$\mathbf{G}$};
\end{tikzpicture}
\end{center}
The edges in this ``tadpole'' have length $1$, except for  $|G C|$  and $|AD|$,
which have lengths $L_{11,2}$ and $L_{11,3}$, respectively.
The angle $\angle CGB = \theta$, and the angles $\angle G B A = \angle B A F$ are $4 \theta$.
The area of the tadpole is $\frac{5}{2}\sin 2\theta + \sin 4\theta$.
The dissection of Fig.~\ref{Fig11S1} and the remaining steps
 in the formation of Fig.~\Ref{Fig11E} can all be carried out with straightedge and compass (assuming,
 of course, that we are given an $11$-gon to start with) .

\begin{figure}[!ht]
\centerline{\includegraphics[angle=270, width=5.5in]{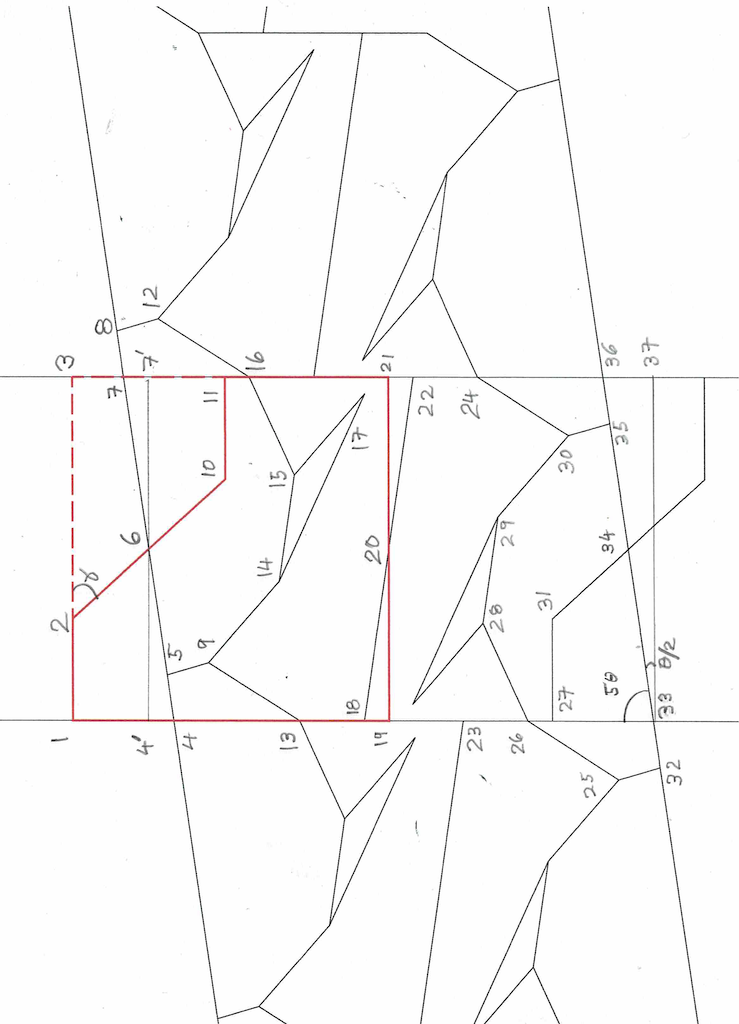}}
\caption{First superposition used to produce the $9$-piece rectangle dissection of Fig.~\ref{Fig11E}.}
\label{Fig11D}
\end{figure} 

\begin{figure}[!ht]
\centerline{\includegraphics[clip=true, trim={0cm, 0cm, 0cm, 0cm},  angle=270,  width=0.8\linewidth]{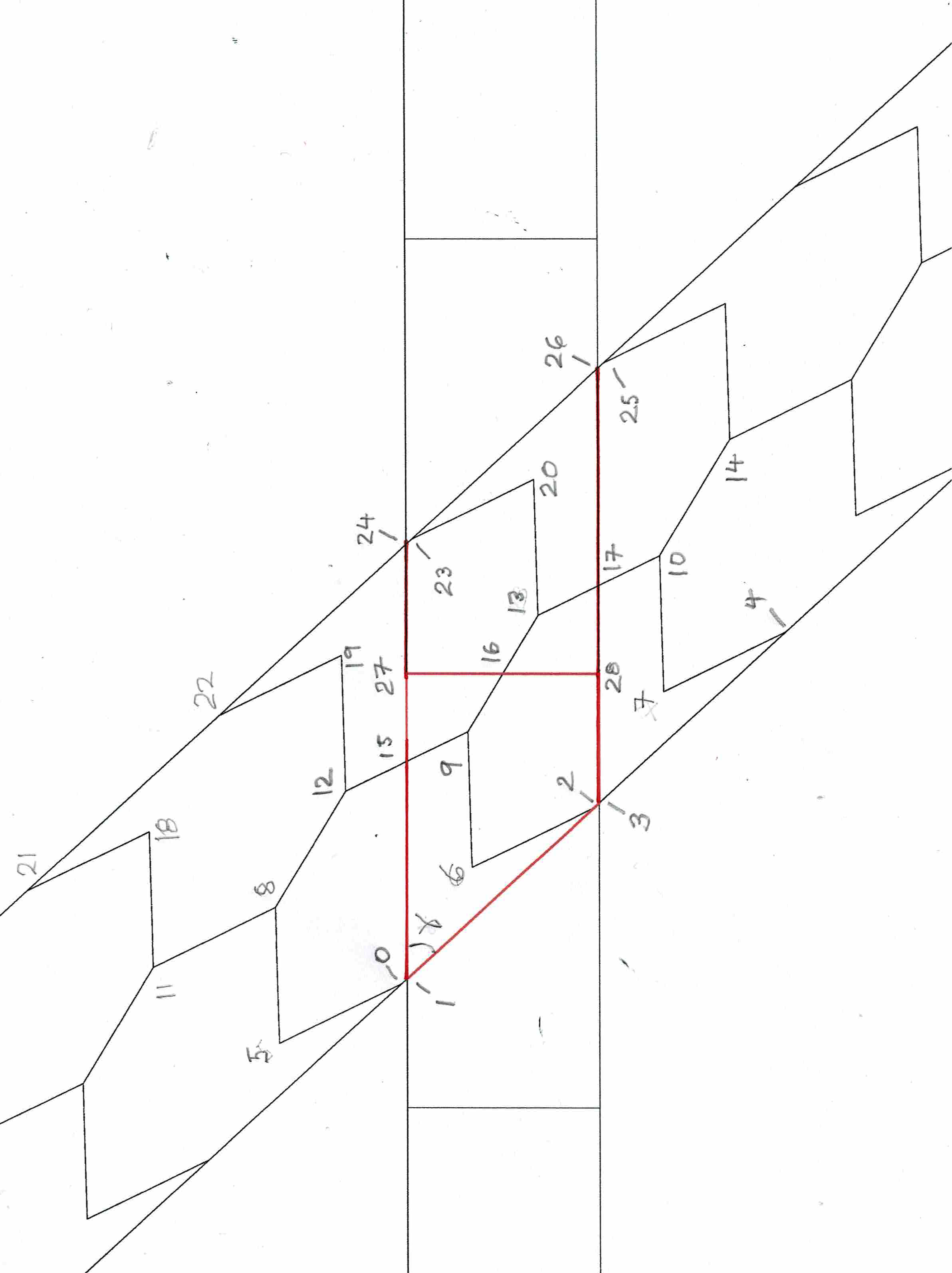}}
\caption{Second superposition used to produce the $9$-piece rectangle dissection of Fig.~\ref{Fig11E}.}
\label{Fig11Db}
\end{figure} 

\begin{figure}[!ht]
\centerline{\includegraphics[clip=true, trim={4cm, 0cm, 4cm, 0cm},  angle=90,  width=0.75\linewidth] 
{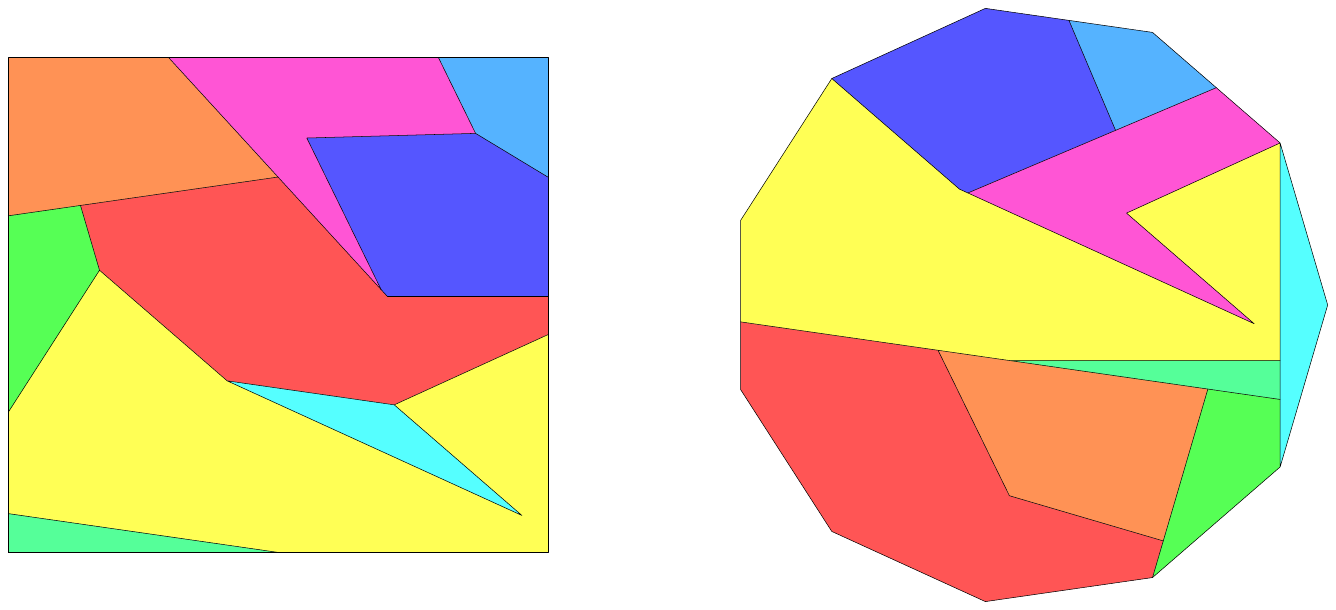}}
\caption{Nine-piece dissection of an  $11$-gon into a rectangle, obtained by combining superpositions in Figs.~\ref{Fig11D} and \ref{Fig11Db}.}
\label{Fig11E}
\end{figure}

The wide, almost horizontal, strip in the first superposition (Fig.~\ref{Fig11D}) is actually a double strip.
Copies of the large heptagonal piece in Fig.~\ref{Fig11S1} are placed along both the bottom edge of the strip
(for example, $35, 30, 29, 28, 26, 25, 32$) and the top edge ($5, 9, 14, 15, 16,  12, 8$).
The interior of the strip is filled with copies of two other pieces from Fig.~\ref{Fig11S1}.
From Fig.~\ref{Fig11S1} we see that $|32, 35| = |18, 22| = L_{11,4}$.

We then superimpose a vertical strip, bounded by the lines $26,13$ 
and $24, 16$. The angle $\angle 26, 33, 36$ between the two strips is $5 \theta$,
and $\angle 36, 33, 37 = \pi/2 - \theta = \theta/2 = \angle 18, 20, 19$.
Along the vertical line $26,13$, the segment $33, 26$ is the side of a 
quadrilateral $33, 26, 25, 32$ with internal angles $2\theta$, $8\theta$, $6\theta$, 
and $6 \theta$, so $|33, 26|  = 4 \cos^2 \theta  - 2 \cos \theta - \frac{3}{5}$. 
From Fig.~\ref{Fig11S1}, $|26,18| =  L_{11,2} - (1-3/5)$, and we know $|18,13| = 3/5$.
Adding up the lengths of the segments, we get  $|33,4| = 8 \cos^2 \theta - 2 \cos \theta -1$.

The vertical strip has width $|19,21| = L_{11,4} \cos(\theta/2) = d_1$ (say) $= 3.1958\ldots$.
This is the width of the final rectangle on the left of Fig.~\ref{Fig11E}.
The height $|1, 19|$  of this rectangle is then determined by the area of the $11$-gon.
This height is $11 \cot \theta/(4 d_1) = d_2$ (say) $= 2.9305\ldots$. 
We easily determine the positions of the points $6, 4^{\prime},
7^{\prime}, 20, 19$, and $21$. The lengths $|2, 10|$ and $|3, 11|$ will be obtained 
from the second 
superposition.\footnote{We need the second superposition to find the angle $\gamma$.}  
We now have full information about the coordinates in Fig.~\ref{Fig11D}.

The second superposition (Fig.~\ref{Fig11Db}) contains a diagonal double strip, formed from copies 
of the ``tadpole'', with a horizontal strip superimposed on it.
We start by constructing the strip of tadpoles, and construct the horizontal strip from it. 
The points are labeled as in Fig.~\ref{Fig11Db}.
The diagonal strip has period $2 \cos \theta$ along the strip (e.g. $|21, 22|$),
and width $3 \sin \theta + 2 \sin 3\theta$ (as can be easily found from the properties of the tadpole)
in the perpendicular direction. The product $A_2$ (say) of these two quantities is the area of a fundamental region for the double strip.

The points $15$, $16$, and $17$ are the midpoints of sides of
the tadpoles. With $15$ as center, we draw a circle of radius $d_1/2$ (taken
from the first superposition), which meets the diagonal
strip at points $1$ and $24$, The line $3, 26$ is constructed similarly.
The trapezoid $1, 3, 28, 27$ will replace the region $2, 10, 11, 3$ in the first superposition.
Let $\gamma$ denote the angle $\angle 15, 1, 3$. The area of this trapezoid can now be 
found in two ways: it is half of $A_2$, that is, $\cos \theta (3 \sin \theta + 2 \sin 3 \theta)$,
and it is also $d_1 \cos \theta \sin \gamma$.
It is also (area of $11$-gon) $-~ d_1 |6,20|$. After some simplification, we find that
$$
\sin \gamma ~=~ \frac{2 \cos^2 \theta - 10 \cos \theta + 11}{2 \cos^2 \theta + 3 \cos \theta  - 4}
~=~ 0.7374\ldots \,.
$$
We can now  give the two sides of the trapezoid that 
were needed for the first superposition: they are $|1,3| = 2 \cos \theta$ and 
$|27,28| = |1,3| \sin \gamma$.

 We now also have full information
about the coordinates in Fig.~\ref{Fig11Db}.
For example, by following around the boundary of the region
$0, 5, 8, 12, 15, 1$, we find that the very short side $0, 1$ of that 
hexagon has 
length $\frac{1}{2}  ( \cos \theta + 2 \cos 3 \theta - d_1 \cos \gamma)$ $=~0.05532\ldots$.
This is the (dark blue) hexagon at the top of the dissected $11$-gon in Fig.~\ref{Fig11E}.

\begin{figure}[!ht]
\centerline{\includegraphics[clip=true, trim={0cm, 4cm, 0cm, 4cm},  angle=0,  width=0.670\linewidth] 
{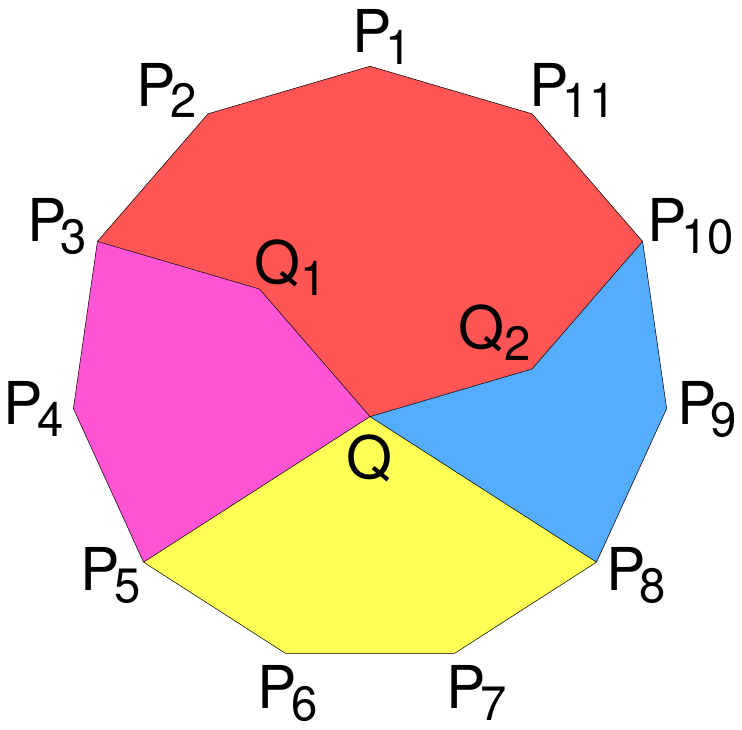}}
\caption{Four-piece dissection of an  $11$-gon used for tessellation in Fig.~\ref{Fig11.2.B}.}
\label{Fig11.2.A}
\end{figure}

\begin{figure}[!ht]
\centerline{\includegraphics[clip=true, trim={0cm, 0cm, 0cm, 0cm},  angle=0,  width=0.60\linewidth] 
{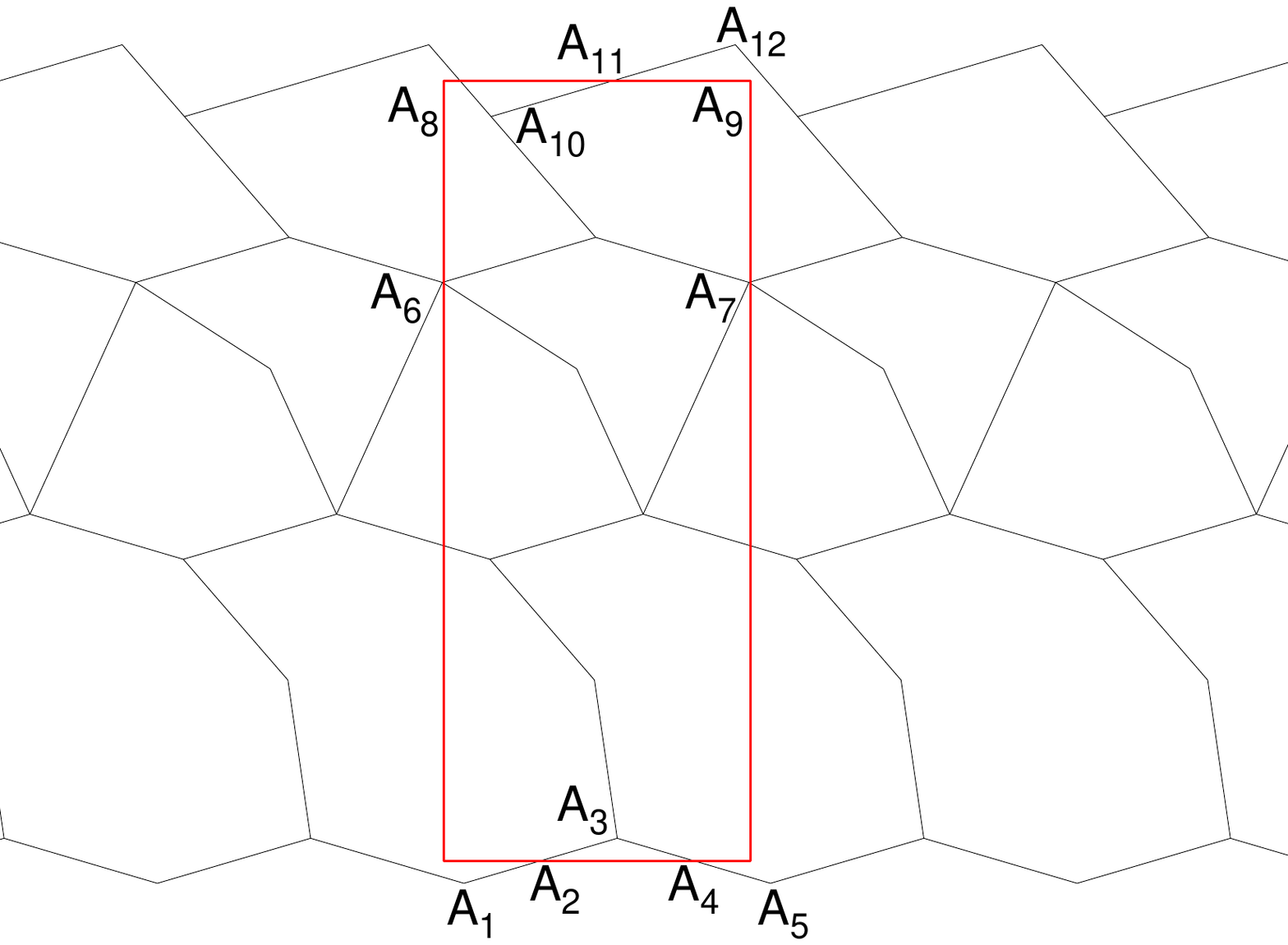}}
\caption{Tessellation of plane built from pieces from Fig.~\ref{Fig11.2.A}.}
\label{Fig11.2.B}
\end{figure}

\subsection{A second nine-piece dissection of an 11-gon to a rectangle}\label{Sec11r2}
Our second nine-piece $11$-gon to rectangle
dissection is slightly simpler, as it only uses a single superposition.
We describe the construction, and give some of the distances and angles,
but leave the detailed verification  of its correctness to the reader.
The starting point is the simple four-piece dissection of the $11$-gon
shown in Fig.~\ref{Fig11.2.A}.
We draw chords from $P_3$ to $P_8$ and from $P_5$ to $P_{10}$,
intersecting at $Q$, say.  $Q$ is located at a distance
$\sin(2 \theta)/ (\cos(\theta) + \cos(2 \theta)) = 0.3002\ldots$ below the center of the $11$-gon.
The segments $P_3 Q$ and $P_{10} Q$ have length $2 \cos \theta$. We replace $P_3 Q$ by a pair of line segments of length $1$, $P_3 Q_1$ and $Q_1 Q$, where $\angle Q_1 P_3 Q = \angle Q_1 Q P_3 = \theta$, 
with a similar construction for $Q_2$.

The angles in Fig.~\ref{Fig11.2.A} are remarkably nice, they are all multiples of $\theta$: 
$\angle Q_1 Q P_5 =  5 \theta$,
$\angle P_5 Q P_8 = 7 \theta$, $\angle P_8 Q Q_2 = 3 \theta$, and so on. 
$Q$ seems to be an auspicious interior point in the $11$-gon.

\begin{figure}[!ht]
\centerline{\includegraphics[clip=true, trim={2cm, 0cm, 2cm, 0cm},  angle=90,  width=0.70\linewidth]
{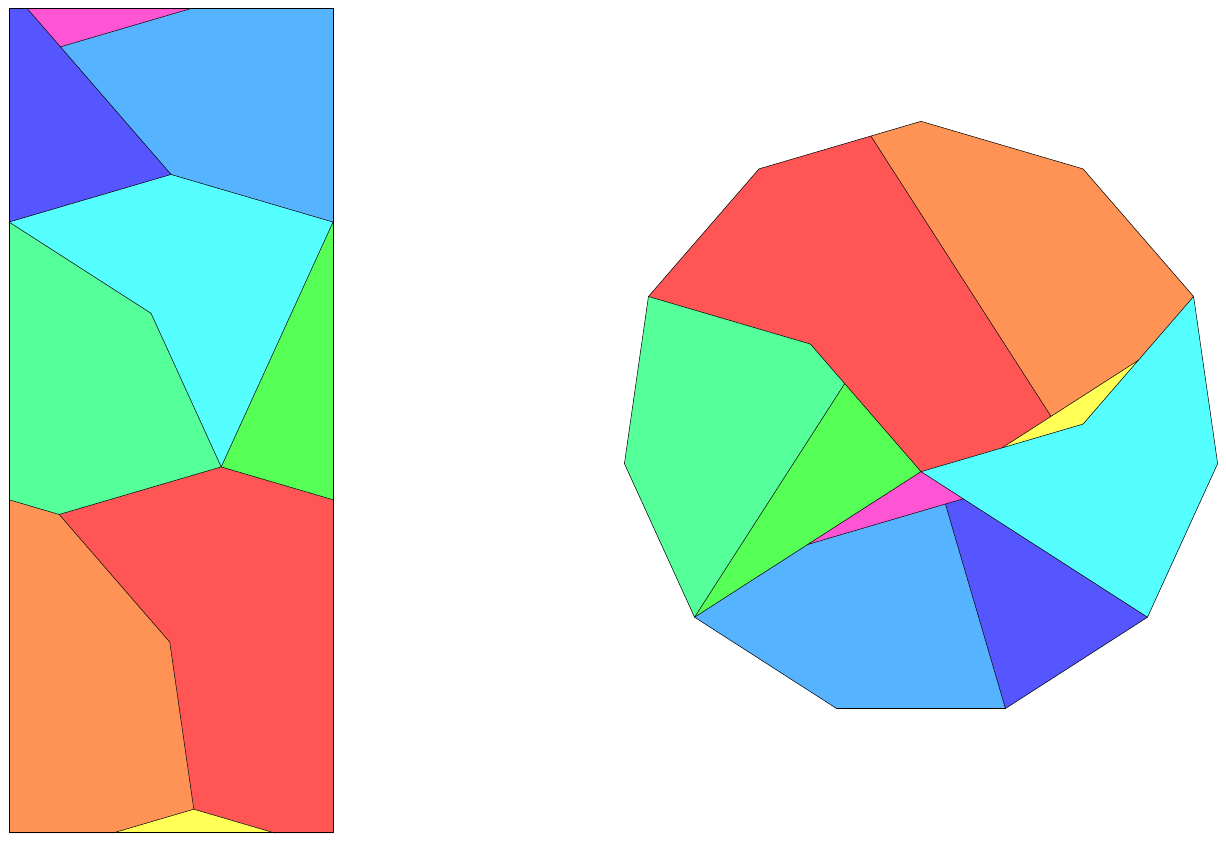}}
\caption{A second nine-piece dissection of an $11$-gon into a rectangle.}
\label{Fig11.2.C}
\end{figure}

We now use these four pieces to build a tessellation of the plane, as shown in Fig.~\ref{Fig11.2.B},
where some of the points have been labeled.
We then cut out the rectangle outlined in red from the tessellation.
The vertical edges of this rectangle  pass through the points $A_6$ and $A_7$, 
so the width of the rectangle (see Fig.~\ref{Fig11.2.A}) is $2 \cos \theta$.
The height is therefore $11/(8 \sin \theta)$.
The rectangle is bounded at the top by a horizontal line through $A_{11}$, the midpoint 
of $A_{10}, A_{12}$, and at the bottom by a
line through the midpoint $A_2$ of $A_1 - A_3$, and 
the midpoint $A_4$ of $A_3 - A_5$. 
This is the rectangle on the left of 
Fig.~\ref{Fig11.2.C}. Finally, the nine pieces in the rectangle can be 
rearranged to form an $11$-gon, as shown on the right in Fig.~\ref{Fig11.2.C}.

\section{Five-piece dissections of a 12-gon}\label{Sec12}

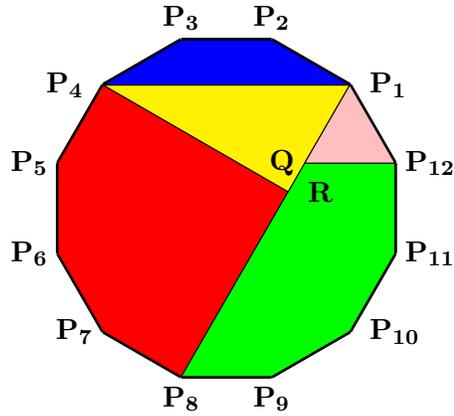
\begin{figure}[!ht]
\begin{center}
\begin{tikzpicture}[scale=1.2]
\coordinate(P2) at  (.5,1.866);
\coordinate(P3) at  (-.5,1.866);
\coordinate(P9) at  (.5,-1.866);
\coordinate(P8) at  (-.5,-1.866);
\coordinate(P1) at (1.366,1.366);
\coordinate(P4) at (-1.366,1.366);
\coordinate(P10) at (1.366,-1.366);
\coordinate(P7) at (-1.366,-1.366);
\coordinate(P12) at (1.866, .5);    
\coordinate(P5) at (-1.866, .5);    
\coordinate(P11) at (1.866, -.5);    
\coordinate(P6) at (-1.866, -.5);    
\coordinate(Q) at (0.866, .5);  
\coordinate(R) at (0.683, .183); 
\draw[ultra thick] (P1) -- (P2);
\draw[ultra thick] (P2) -- (P3);
\draw[ultra thick] (P3) -- (P4);
\draw[ultra thick] (P4) -- (P5);
\draw[ultra thick] (P5) -- (P6);
\draw[ultra thick] (P6) -- (P7);
\draw[ultra thick] (P7) -- (P8);
\draw[ultra thick] (P8) -- (P9);
\draw[ultra thick] (P9) -- (P10);
\draw[ultra thick] (P10) -- (P11);
\draw[ultra thick] (P11) -- (P12);
\draw[ultra thick] (P12) -- (P1);

\draw[ultra thick] (P4) -- (P1);
\draw[ultra thick] (P8) -- (P1);
\draw[ultra thick] (Q) -- (P12);
\draw[ultra thick] (P4) -- (R);
\draw[fill = red] (P4) -- (P5) -- (P6) -- (P7) -- (P8) -- (R) -- (P4);
\draw[fill = green] (P8) -- (P9)-- (P10) -- (P11) -- (P12) -- (Q) -- (P8);
\draw[fill = yellow] (P4) -- (R) -- (P1) -- (P4);
\draw[fill = blue] (P4) -- (P1) -- (P2) -- (P3)  -- (P4);
\draw[fill = pink] (P12) -- (P1) -- (Q) -- (P12);

\node[right] at (P1) {$~\mathbf{P_1}$};
\node[above] at (P2) {$\mathbf{P_2}$};
\node[above] at (P3) {$\mathbf{P_3}$};
\node[left] at (P4) {$\mathbf{P_4}~$};
\node[left] at (P5) {$\mathbf{P_5}$};
\node[left] at (P6) {$\mathbf{P_6}$};
\node[left] at (P7) {$\mathbf{P_7}$};
\node[below] at (P8) {$\mathbf{P_8}$};
\node[below] at (P9) {$\mathbf{P_9}$};
\node[right] at (P10) {$~\mathbf{P_{10}}$};
\node[right] at (P11) {$\mathbf{P_{11}}$};
\node[right] at (P12) {$\mathbf{P_{12}}$};
\node[left] at (Q) {$\mathbf{Q}$};
\node[right] at (R) {$~\mathbf{R}$};
\end{tikzpicture}
\caption{A five-piece dissection of a $12$-gon.}
\label{Fig12L}
\end{center}
\end{figure}

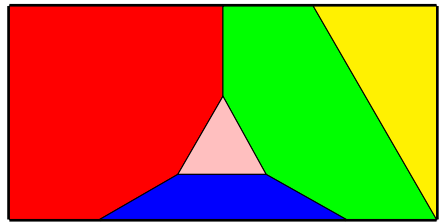
\begin{figure}[!ht]
\begin{center}
\begin{tikzpicture}[scale=1.2]
\def \w{4.732};
\def \h{2.366};
\coordinate(a1) at  (0,0);
\coordinate(a2) at  (0,\h);
\coordinate(a5) at  (\w,\h);
\coordinate(a6) at  (\w,0);
\coordinate(a3) at  (\h,\h);
\coordinate(a8) at  (1,0);
\coordinate(a7) at  (\w - 1,0);
\coordinate(a4) at  (\w - 1.366, \h);
\coordinate(a11) at  (\h ,\h-1);
\coordinate(a9) at  (1.866,.5);
\coordinate(a10) at  (\w - 1.888, .5);

\draw[ultra thick] (a1) -- (a2);
\draw[ultra thick] (a2) -- (a5);
\draw[ultra thick] (a5) -- (a6);
\draw[ultra thick] (a6) -- (a1);
\draw[ultra thick] (a8) -- (a9);
\draw[ultra thick] (a9) -- (a11);
\draw[ultra thick] (a9) -- (a10);
\draw[ultra thick] (a10) -- (a11);
\draw[ultra thick] (a10) -- (a7);
\draw[ultra thick] (a11) -- (a3);
\draw[ultra thick] (a4) -- (a6);

\draw[fill=red] (a1) -- (a2) -- (a3) -- (a11) -- (a9) -- (a8) -- (a1);
\draw[fill=green] (a3) -- (a4) -- (a6) -- (a7) -- (a10) -- (a11) -- (a3);
\draw[fill=pink] (a9) -- (a10) -- (a11) -- (a9);
\draw[fill=yellow] (a4) -- (a5) -- (a6) -- (a4);
\draw[fill=blue] (a7) -- (a8) -- (a9) -- (a10) -- (a7);;
\end{tikzpicture}
\caption{The pieces reassembled to form a rectangle.}
\label{Fig12B}
\end{center}
\end{figure}

\begin{figure}[!ht]
\centerline{\includegraphics[angle=0, width=4.5in]{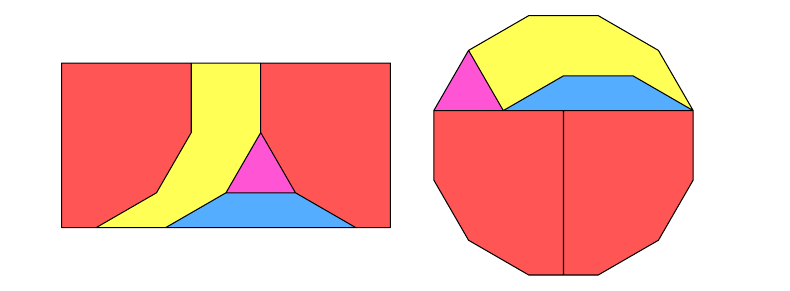}}
\caption{An alternative five-piece dissection of a dodecagon.}
\label{Fig12.2}
\end{figure}

We start with a $12$-gon with edge-length $1$ (see Fig.~\ref{Fig12L}).
Draw chords from $P_1$ to $P_4$ and $P_8$.
Draw a perpendicular from $P_4$  to $P_1 - P_8$, meeting it at $R$, 
and draw an equilateral triangle $P_1 Q P_{12}$
that touches $P_1 - P_8$. 

The angle $\angle P_1 P_4 R$ is $\pi/6 = 30^\circ$.
The lengths of the line segments are as follows: 
$|P_1 P_4| = |Q P_8| = L_{3,12} = 1 + \sqrt{3}$,
$|P_4 R| = |R P_8| = (3 + \sqrt{3})/2$,
and $|P_1 R| = (1+\sqrt{3})/2$.

After the pieces are rearranged (Fig.~\ref{Fig12B}),
the resulting rectangle has width $w = 3 + \sqrt{3}$
and height $h= (3 + \sqrt{3})/2$. We then easily check that the product $w h$ is equal 
to the area $3 \cot 15^\circ$.

A second dissection of the $12$-gon (although with a non-convex piece) is 
shown in Fig.~\ref{Fig12.2}.

\section{A seven-piece dissection of a 14-gon 
and a nine-piece dissection of a 16-gon}\label{Sec14}


\begin{figure}[!ht]
\centerline{\includegraphics[angle=0, width=3.0in]{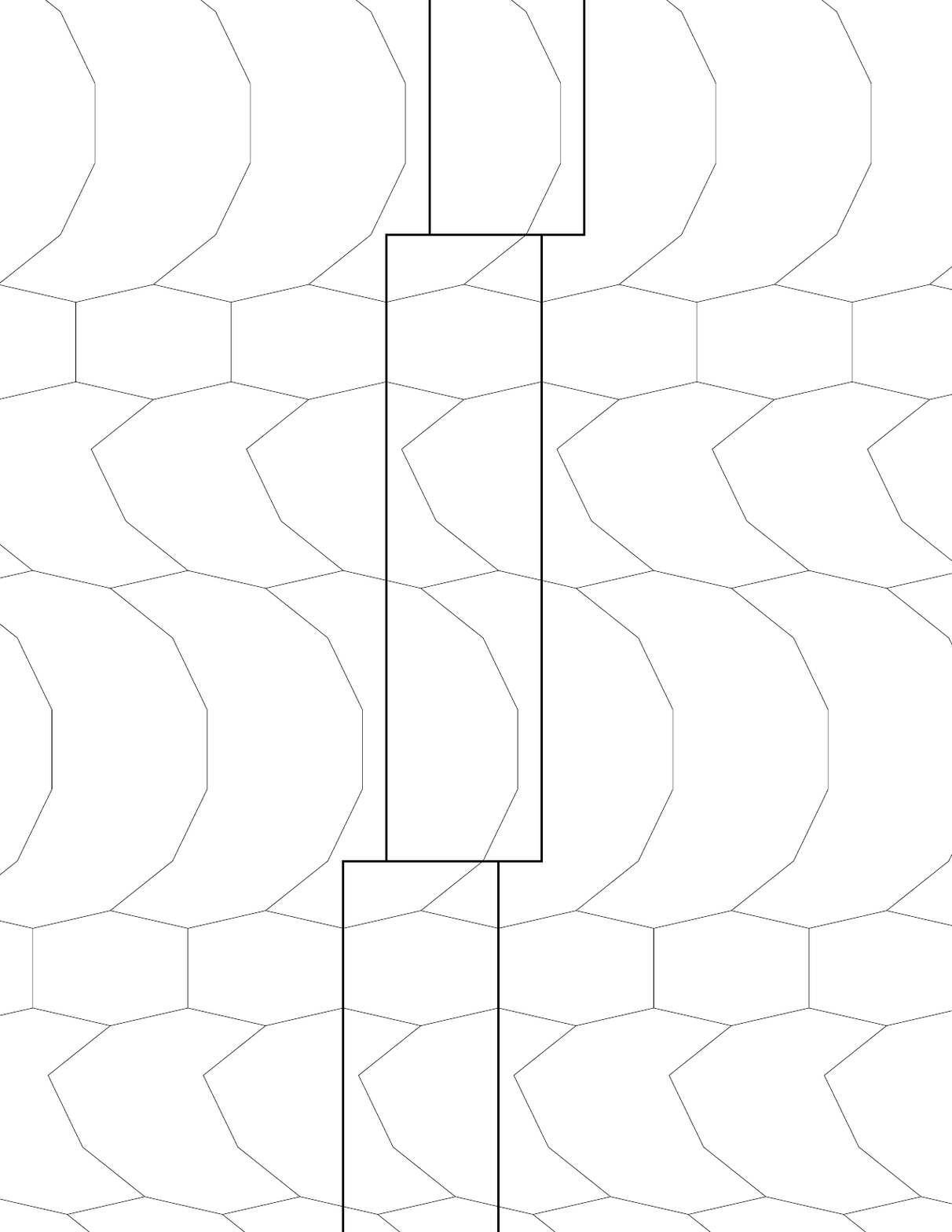}}
\caption{A tessellation based on a $14$-gon, with the rectangle that leads to the dissection in Fig.~\ref{Fig14B}.}
\label{Fig14A}
\end{figure}

\begin{figure}[!ht]
\centerline{\includegraphics[angle=0, width=4.0in]{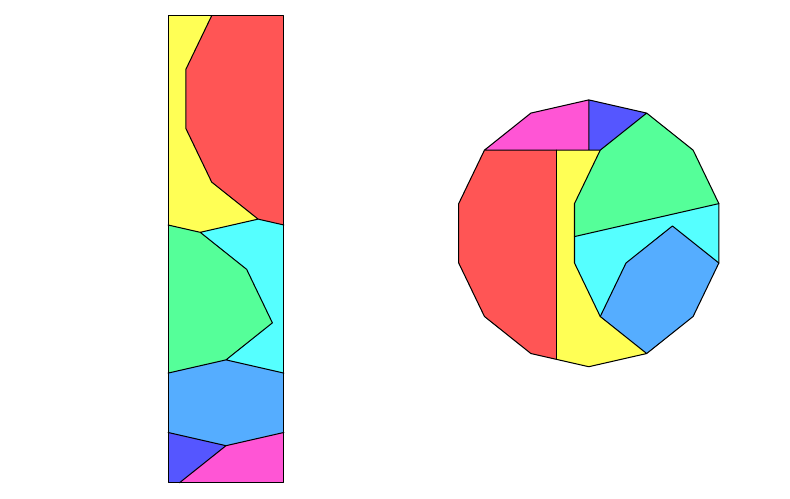}}
\caption{A seven-piece dissection of a $14$-gon into a rectangle.} 
\label{Fig14B}
\end{figure}

\begin{figure}[!ht]
\centerline{\includegraphics[angle=0, width=3.0in]{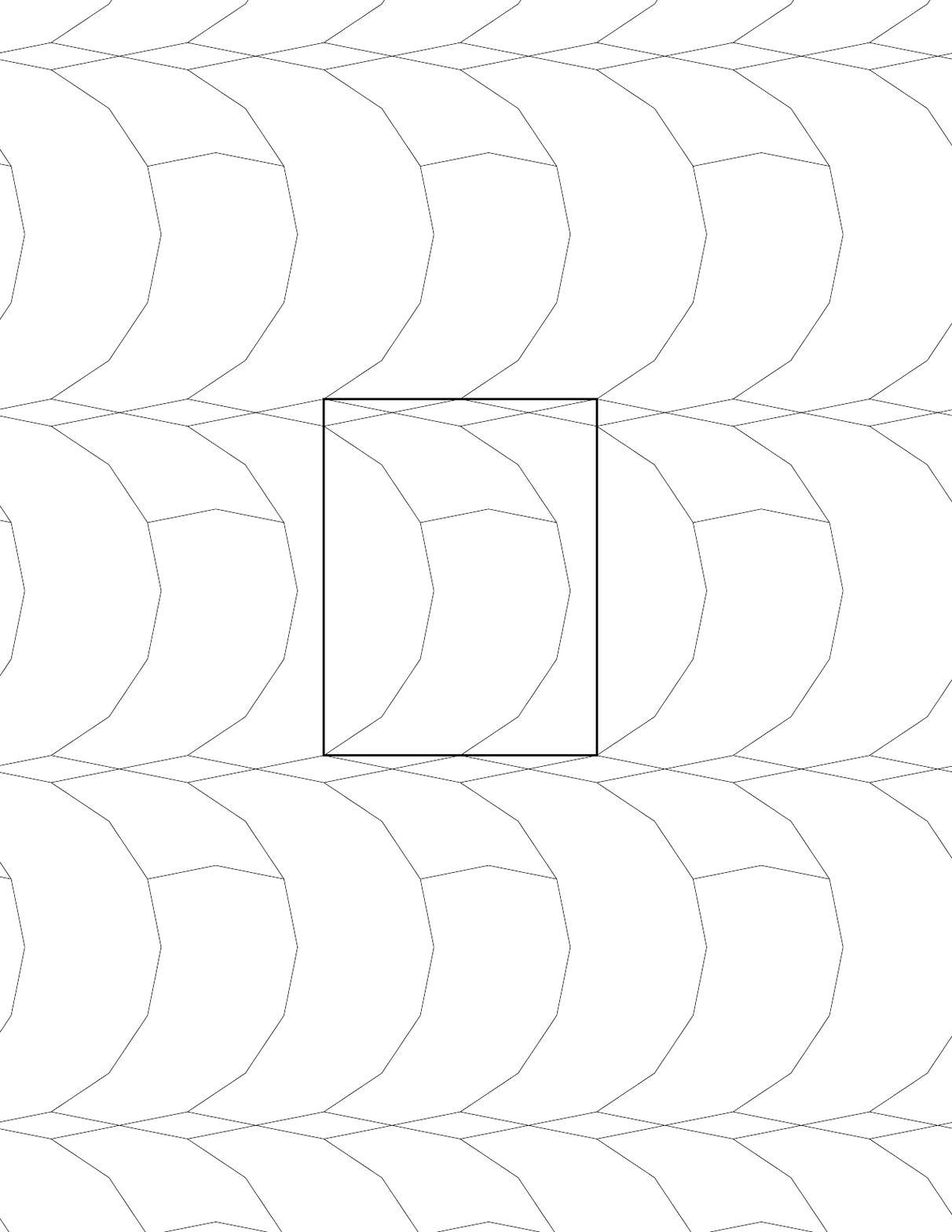}}
\caption{A tessellation based on a $16$-gon, with the rectangle that leads to the dissection in Fig.~\ref{Fig16B}.}
\label{Fig16A}
\end{figure}

\begin{figure}[!ht]
\centerline{\includegraphics[angle=0, width=4.0in]{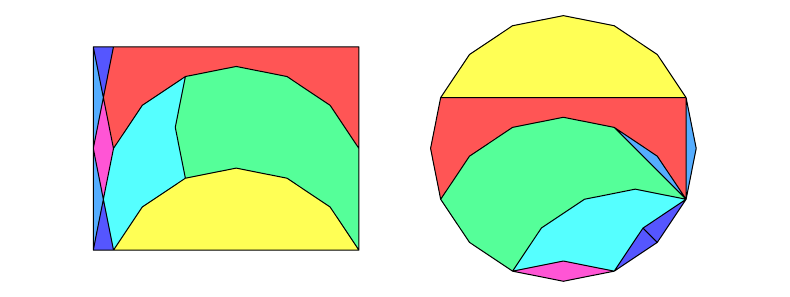}}
\caption{A nine-piece dissection of a $16$-gon into a rectangle.} 
\label{Fig16B}
\end{figure}

\begin{figure}[htb]
        \begin{minipage}[b]{0.45\linewidth} 
        \centering
        \includegraphics[width=0.8\textwidth]{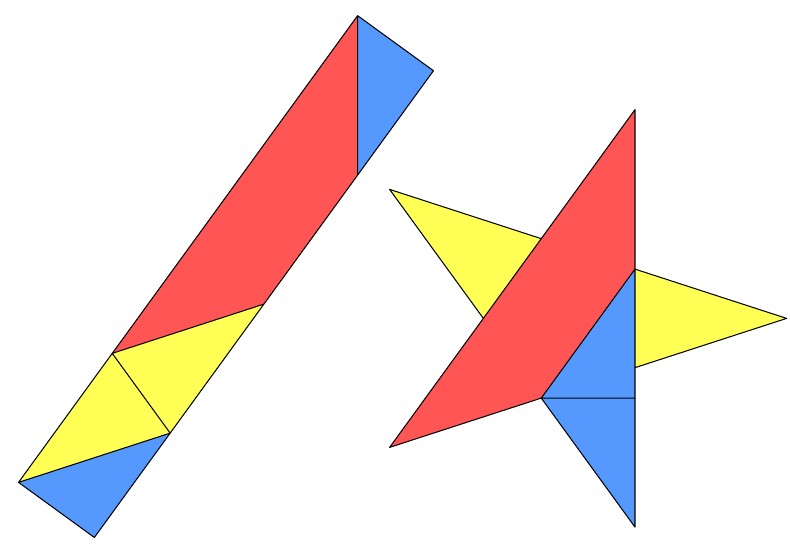}
        \caption{A five-piece dissection of a $\{5/2\}$ pentagram.}
        \label{FigPentagram}
        \end{minipage}
        \hspace{0.5cm}
        \begin{minipage}[b]{0.45\linewidth}
        \centering 
        \includegraphics[width=0.8\textwidth]{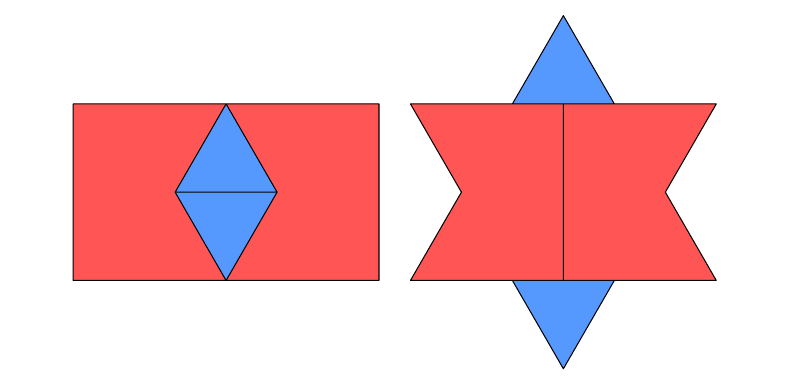}
        \caption{A four-piece dissection of a $\{6/2\}$ hexagram.}
        \label{FigHexagram}
        \end{minipage}
        \end{figure}

It is known that $s(14) \le 9$ and $s(16) \le 10$ \cite{GDDb}. Figure~\ref{Fig14A} shows a tessellation of the plane 
based on a $14$-gon, and a rectangle superimposed on it
which leads to the seven-piece dissection shown in Fig.~\ref{Fig14B}.
Figures~\ref{Fig16A} and \ref{Fig16B} play  similar roles for the $16$-gon.

These will be our final examples of a regular polygon to rectangle dissections.
Other examples with larger numbers of sides may be found in
the {\em Rectangle Dissections} section of \cite{GDDb}, and in a projected sequel to the present work.

\section{Selected dissections of star polygons to rectangles}\label{SecStar}

We give four examples of especially elegant dissections of star polygons to rectangles.
These are taken from \cite{GDDb}, where many further examples can be found.


\begin{figure}[htb]
        \begin{minipage}[b]{0.45\linewidth} 
        \centering
        \includegraphics[width=0.8\textwidth]{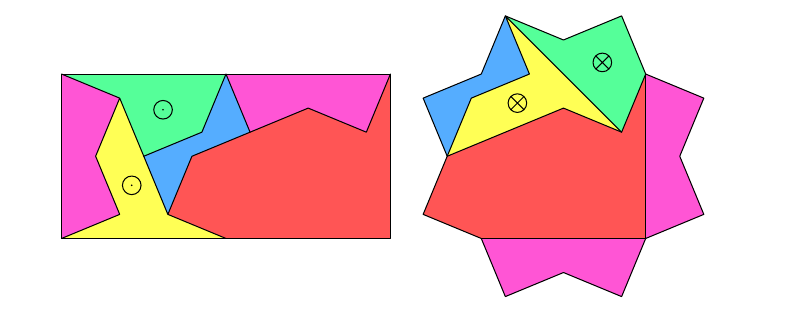}
        \caption{A six-piece dissection of an $\{8/2\}$ octagram. Two pieces must be turned over.}
        \label{FigOctagram82}
         \end{minipage}
        \hspace{0.5cm}
        \begin{minipage}[b]{0.45\linewidth}
        \centering 
        \includegraphics[width=0.8\textwidth]{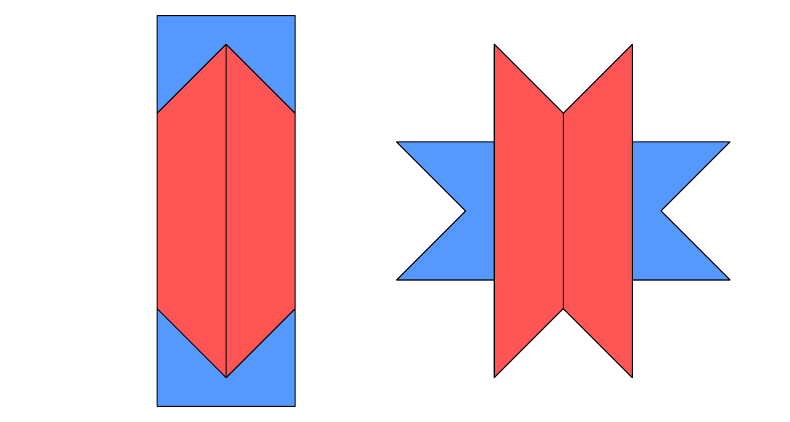}
        \caption{A four-piece dissection of an $\{8/3\}$ octagram.}
        \label{FigOctagram83}
        \end{minipage}
        \end{figure}

\section{Three-piece dissections of a Greek cross}\label{SecGC}

Most authors who  study dissections of polygons include the Greek cross,
so we briefly discuss it here. The classical four-piece dissection of a Greek cross 
into a square can be seen for example in~\cite[Fig.~9.1]{Lind64}. 
The dissected square has $4$-fold rotational symmetry.

\begin{figure}[htb]
        \begin{minipage}[b]{0.45\linewidth} 
        \centering
        \includegraphics[width=0.8\textwidth]{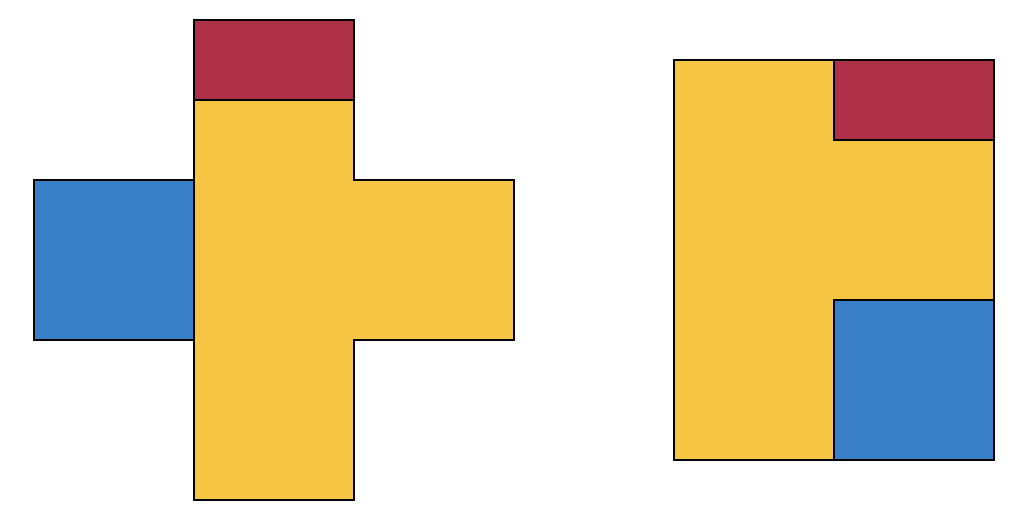}
        \caption{Eppstein's $3$-piece dissection of a Greek cross into a rectangle \cite{Epp22}.}
        \label{FigEpp1}
        \end{minipage}
        \hspace{0.5cm}
        \begin{minipage}[b]{0.45\linewidth}
        \centering
        \includegraphics[angle=180, width=0.9\textwidth]{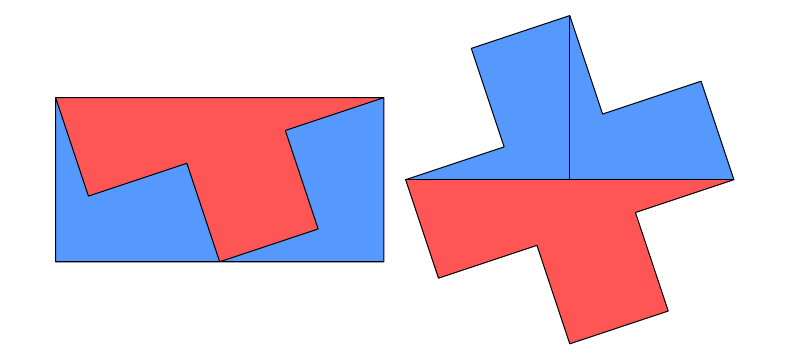}
        \caption{A $3$-piece dissection of a Greek cross into a rectangle  (from \cite{GDDb}). }
        \label{FigTh11}
        \end{minipage}
        \end{figure}

Three pieces seem to be the minimal number needed to form a
rectangle from a Greek cross. The simplest three-piece construction cuts off two
opposite arms from the cross  and places them at the ends of the other two arms, forming 
a $1 \times 5$ rectangle
\begin{tikzpicture}[scale=.35]
\draw (0,0)--(5,0); 
\draw (0,1)--(5,1);
\draw (0,0)--(0,1);
\draw (1,0)--(1,1);
\draw[dashed] (2,0)--(2,1);
\draw[dashed] (3,0)--(3,1);
\draw (4,0)--(4,1);
\draw (5,0)--(5,1);
\end{tikzpicture}\,.
~~Eppstein~\cite{Epp22} gives a  three-piece dissection into non-convex pieces,
shown in Fig.~\ref{FigEpp1}, and the database \cite{GDDb} gives another  (Fig.~\ref{FigTh11}),
similar in spirit to the four-piece dissection into a square.


\section{Curved cuts are sometimes essential}\label{SecCurved}
We know of no theorem  which  will guarantee that  polygonal cuts are sufficient
to achieve $s(n)$ or $r(n)$.
The following are three examples of other situations  where it seems clear that minimal dissections
can {\em not} be achieved using only polygonal cuts.

\begin{figure}[!ht]
\centerline{\includegraphics[angle=0, width=2.0in]{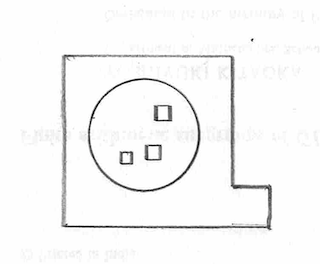}}
\caption{A two-piece dissection that can be accomplished by a single circular cut, but surely not by a polygonal cut. The polygon contains  three small squares holes
(Richard C. Schroeppel and Andy Latto).}
\label{FigRCS}
\end{figure} 

1. Take a square with a smaller square attached to it, and
cut out three small square holes at random positions in the interior.  Call this figure $A$.
For figure $B$, make a circular cut enclosing the three holes, and rotate the interior of the circle by
a small  random angle.  This gives a two-piece dissection of $A$ to $B$ which surely cannot be
accomplished with a single  polygonal cut.  This example was suggested by 
Richard C. Schroeppel and Andy Latto, (personal communication).

\begin{figure}[!ht]
\centerline{\includegraphics[angle=0, width=3.7in]{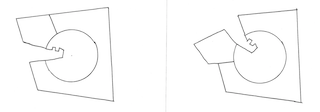}}
\caption{A three-piece dissection that appears to require a curved piece (David desJardins).}
\label{FigDDJ}
\end{figure}

2. Figure~\ref{FigDDJ} shows an example due to David desJardins (personal communication)
of a three-piece dissection between two simply connected polygonal regions
that appears to require a curved piece.

\begin{figure}[!ht]
\centerline{\includegraphics[angle=0, width=5in]{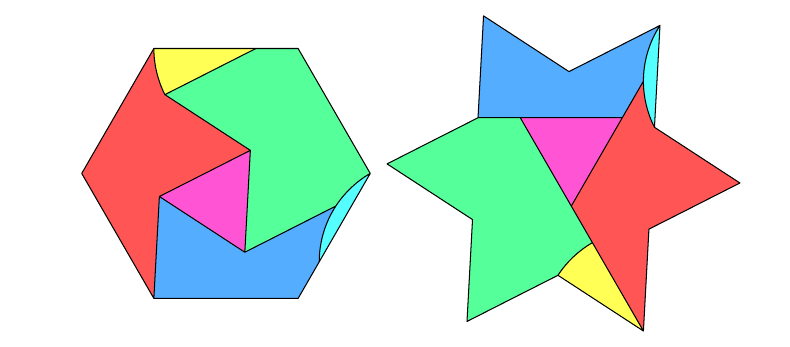}}
\caption{A $6$-piece hexagon to hexagram dissection that avoids turning pieces
over, but uses a curved cut.}
\label{FigHexHex}
\end{figure}

3.  Greg Frederickson \cite[Fig.\ 13.6]{Fred97} gives an example of $6$-piece dissection
of a hexagon into a hexagram which requires that two of the pieces
be turned over. This can be modified to avoid turning the pieces over at the
cost of adding an extra piece. But if  curved cuts are used, this can be accomplished without adding
the extra piece, as shown in Fig.~\ref{FigHexHex}.
We conjecture that a $6$-piece hexagon to hexagram dissection that avoids turning pieces over cannot be
constructed using only polygonal cuts.

\begin{figure}[htb]
        \begin{minipage}[b]{0.45\linewidth} 
        \centering
        \centerline{\includegraphics[clip=true, trim={0cm, 3cm, 0cm, 3cm},  angle=0,  width=0.8\linewidth]{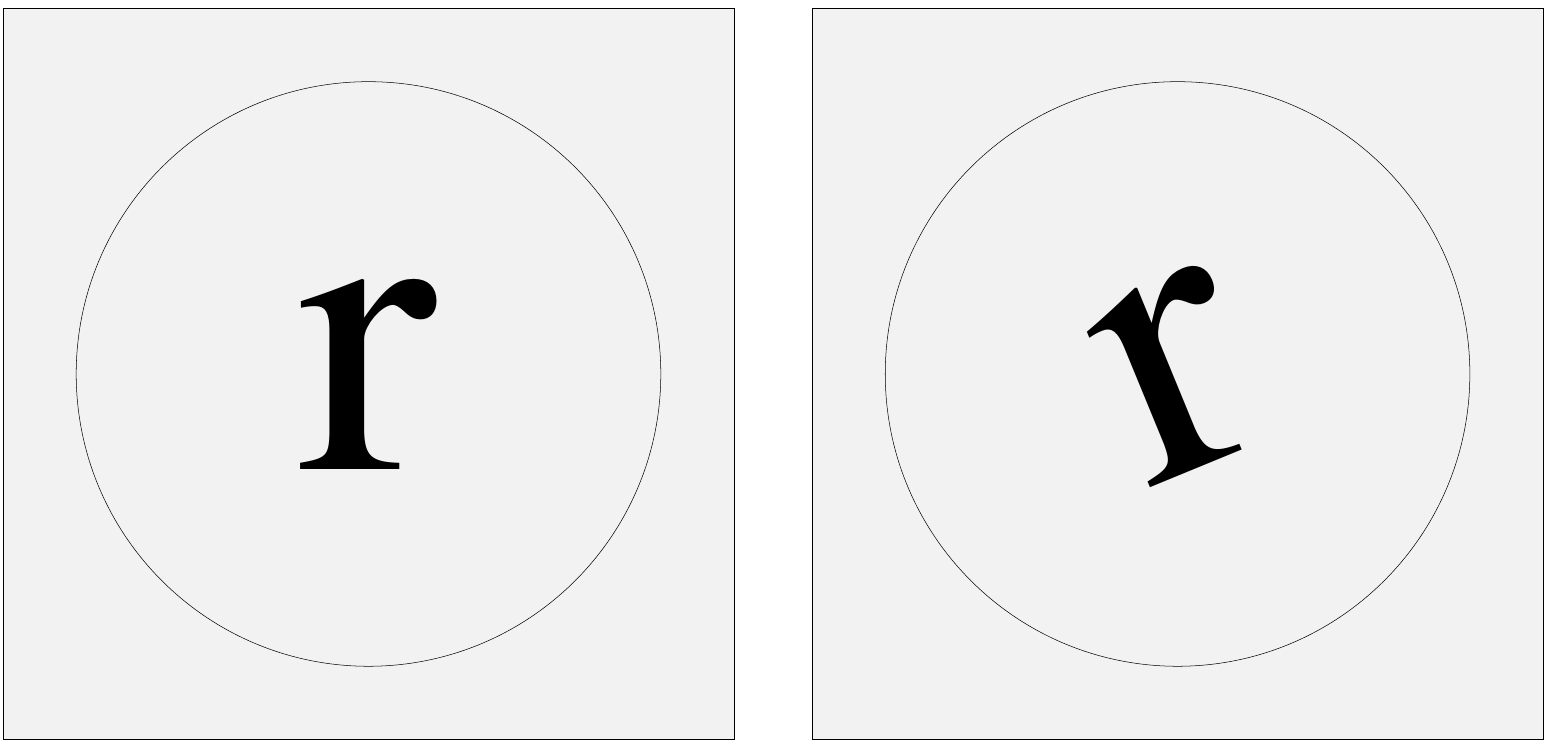}}  
       \caption{A rotation produced by a single circular cut ... }
        \label{FigRota}
        \end{minipage}
        \hspace{0.5cm}
        \begin{minipage}[b]{0.45\linewidth}
        \centerline{\includegraphics[clip=true, trim={0cm, 3cm, 0cm, 3cm},  angle=0,  width=0.8\linewidth]{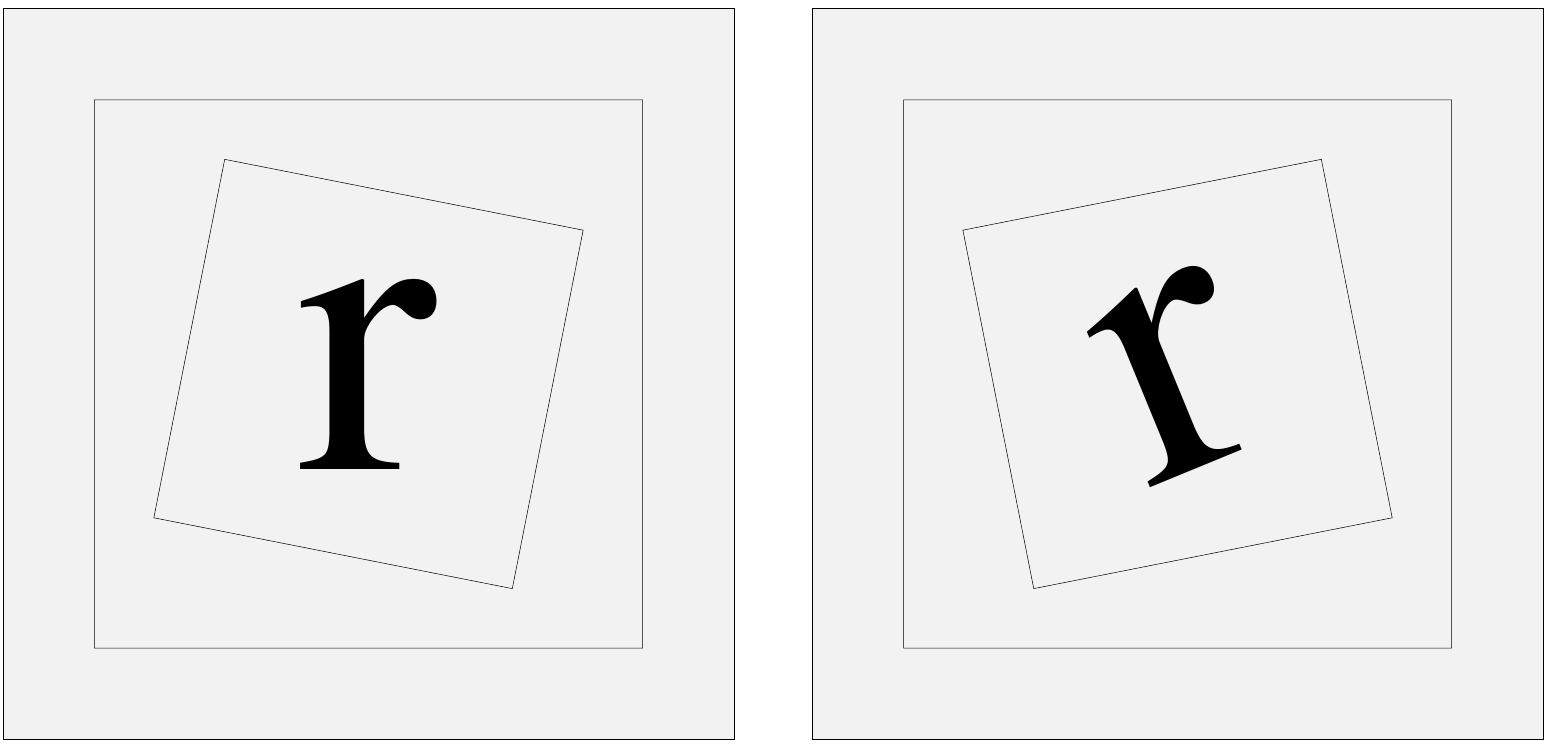}}  
        \caption{... can also be achieved with two square cuts and turning a piece over.}
        \label{FigRotd}
        \end{minipage}
        \end{figure}

4.  In this regard, it is worth pointing out that a rotation by any desired angle that uses a 
single circular cut (see Fig.~\ref{FigRota}) can be accomplished by two square cuts and 
turning a   piece over  (Fig.~\ref{FigRotd}).   


\section{Acknowledgments}
Thanks to Adam Gsellman for telling us about his polygon to
rectangle dissections, which was the seed that led to the present paper.
Thanks also  to  David desJardins, Andy Latto, Richard C. Schroeppel, and Allan C. Wechsler for helpful comments.
The writing of this paper depended heavily on PostScript, LaTeX, Tikz, Maple, WolframAlpha, and email.

\bigskip
\hrule
\bigskip

\noindent 2020 Mathematics Subject Classification 52B45

\end{document}